\documentclass[12pt]{article}

\usepackage{amsfonts}
\usepackage{amsmath}
\usepackage{amssymb}
\usepackage{graphicx}
\usepackage[usenames]{color}

\usepackage{mathrsfs}

\newtheorem{theorem}{Theorem}[section]
\newtheorem{proposition}[theorem]{Proposition}

\newtheorem{lemma}[theorem]{Lemma}
\newtheorem{remark}[theorem]{Remark}
\newtheorem{definition}[theorem]{Definition}

\numberwithin{equation}{section} 

\def\Sc{{Schr\"odinger} }

\def\en t{{{\rm Z}\mkern-5.5mu{\rm Z}}}

\def\oo{{\omega}}

\def\dd{\displaystyle}

\def\<{\left<}
\def\>{\right>}
\def\({\left(}
\def\){\right)}

\def\9{{\infty}}
\def\barr{\begin{array}}
\def\earr{\end{array}}

\def\wt{\widetilde}
\def\wh{\widehat}

\def\lbb{{\lambda}}
\def\g{{\gamma}}
\def\a{{\alpha}}

\def\vsp{\vspace*{1,5mm}\\ }

\def\n{\noindent }

\def\dd{\displaystyle}

\def\3{\subset }
\def\na{{\nabla}}
\def\sk{\smallskip }

\def\bk{\bigskip }
\def\e{{\epsilon}}
\def\ve{{\varepsilon}}

\begin{document}

\begin{center}
{\Large{\bf Stochastic nonlinear \Sc equations}}
\bigskip\bk

{\large{\bf Viorel Barbu}}\footnote{Octav Mayer Institute of
Mathematics (Romanian Academy)   and Al.I. Cuza University and,
700506, Ia\c si, Romania.  This author was supported by a grant of
the Romanian National Authority for Scientific Research,
PN-II-ID-PCE-2011-3-0027 and BiBos - Research Centre.}, {\large{\bf
Michael R\"ockner}}\footnote{Fakult\"at f\"ur Mathematik,
Universit\"at Bielefeld,  D-33501 Bielefeld, Germany. This research
was supported by the DFG through CRC 701.}, {\large{\bf Deng
Zhang}}\footnote{Fakult\"at f\"ur Mathematik, Universit\"at
Bielefeld,  D-33501 Bielefeld, Germany,  and Insti\-tute of Applied
Mathematics, Academy of Mathematics and Systems Science, Chinese
Academy of Sciences. Beijing 100190, China.}
\end{center}

\bk\bk\bk

\begin{quote}
\n{\small{\bf Abstract.} This paper is devoted to the well-posedness
of stochastic nonlinear Schr\"odinger equations in the energy space
$H^1(\mathbb{R}^d)$, which is a natural continuation of our recent
work \cite{BRZ14}. We consider both focusing and defocusing
nonlinearities and prove global well-posedness in
$H^1(\mathbb{R}^d)$, including also the pathwise continuous
dependence on initial conditions, with exponents exactly the same as
in the deterministic case. In particular, this work improves earlier
results in \cite{BD03}. Moreover, the local existence, uniqueness
and blowup alternative are also established for the energy-critical
case.  The approach presented here is mainly based on the rescaling
approach already used in \cite{BRZ14} to study the  $L^2$ case and
also on the Strichartz estimates established in \cite{MMT08} for
large perturbations of
the Laplacian.} \\

\end{quote}

\vfill

\section{Introduction and main results} \label{Intro}

Let us consider the stochastic nonlinear Schr\"{o}dinger equation
with linear multiplicative noise
\begin{align} \label{equax}
      &idX(t,\xi)=\Delta X(t,\xi)dt+\lambda|X(t,\xi)|^{\alpha-1}X(t,\xi)dt \nonumber \\
      &\qquad \qquad \quad -i\mu(\xi) X(t,\xi)dt+iX(t,\xi)dW(t,\xi),\
      t\in(0,T),\ \xi\in \mathbb{R}^d,
      \\
      &X(0)=x.  \nonumber
\end{align}
Here $\lambda=\pm 1$, $\alpha>1$ and $W$ is the colored Wiener
process
\begin{equation} \label{W}
     W(t,\xi)=\sum^N\limits_{j=1}\mu_j e_j(\xi)\beta_j(t),t\geq
     0,\xi\in \mathbb{R}^d,
\end{equation}
with $\mu_j \in \mathbb{C}$, $e_j(\xi)$ real-valued functions and $
\beta_j(t)$ independent real Brownian motions on a probability space
$(\Omega,\mathcal {F},\mathbb{P})$ with natural filtration
$(\mathcal {F}_t)_{t\geq 0}$, $1\leq j \leq N$. In this paper for
simplicity we assume $N<\9$.

As in the physical context \cite{BG09}, we choose $\mu$ of the form
\begin{equation*}
     \mu(\xi)=\frac{1}{2}\sum^N\limits_{j=1}|\mu_j|^2e_j^2(\xi),\xi\in\mathbb{R}^d,
\end{equation*}
so that $|X(t)|_2^2$ is a martingale, from which one can
define the ``physical probability law''(see \cite{BG09}).\\

In the deterministic case $\mu_j=0$, $1\leq j \leq N$, it is well
known (see \cite{K89,LP09}) that \eqref{equax} is globally well
posed in $H^1(\mathbb{R}^d)$ in the defocusing case $\lbb=-1$ with
the subcritical exponents of the nonlinearity
\begin{equation} \label{a-defo}
    1<\a< 1+\frac{4}{(d-2)_+},
\end{equation}
while in the focusing case $\lbb=1$ with the exponents
\begin{align} \label{a-fo}
    1<\a < 1+ \frac{4}{d}.
\end{align}
Here $1+\frac{4}{(d-2)_+}=1+\frac{4}{d-2}\ (resp.\ \9)$ with $d\geq
3\ (resp.\ d=1,2).$

In the stochastic case, the authors in \cite{BD03} (see also
\cite{BD99}) studied the conservative case $Re \mu_j=0$, $1\leq
j\leq N$, i.e. $W$ is a purely imaginary noise. They proved the
local existence and uniqueness with $\alpha$ satisfying
$$\left\{
  \begin{array}{ll}
    1<\alpha< \9, & \hbox{if $d=1$ or $2$;} \\
    1<\alpha<5, & \hbox{if $d=3$;} \\
    2\leq \alpha < 1+\frac{4}{d-2}, & \hbox{if $d=4,5$;}\\
    \a<1+\frac{2}{d-1},& \hbox{if $d\geq6$;}
  \end{array}
\right.$$ and then the global well-posedness under the further
assumptions that $\a<1+\frac{4}{d}$ or $\lambda=-1$. Hence, when $d
\geq 6$, the global well-posedness is established only for the
restrictive exponents $\a < 1+ \frac{2}{d-1}$. We also refer to
\cite{BF12} for stochastic nonlinear Schr\"odinger equation
(however, only for one-dimension
noise) with real-valued potentials in the conservative case. \\

The starting point of this article is our recent work \cite{BRZ14},
where we obtain the global well-posedness of (\ref{equax}) in $L^2$
space with exponent $\alpha\in (1,1+\frac{4}{d})$, i.e. in the same
range as in the deterministic case.

The main aim of the present work is to study the global
well-posedness of (\ref{equax}) in $H^1(\mathbb{R}^d)$ with general
$\mu_j\in\mathbb{C}$ as in the physical context \cite{BG09},
including the non-conservative case. We prove the global
well-posedness, including also the pathwise continuous dependence on
initial conditions, with $\a$ in the ranges (\ref{a-defo}) and
(\ref{a-fo}) in the defocusing and focusing cases respectively, i.e.
in exactly the same ranges as in the deterministic case. In
particular, these sharper results fill the gap for $\a$ in
\cite{BD03} mentioned above.

Moreover, the local well-posedness is also established in Section
\ref{LWP} for the energy-critical case $\lbb=\pm 1$,
$\a=1+\frac{4}{d-2}$ with $d\geq 3$ and also for the focusing
mass-(super)critical case $\lbb= 1$, $1+\frac{4}{d} \leq \a
<1+\frac{4}{(d-2)_+}$ with $d\geq 1$. The local results established
in the latter case allow to study the noise effect on blowup
phenomena,
which will be contained in forthcoming work.\\

Before we show the main global well-posedness result, let us first
present the spatial decay assumption on $\{e_j\}_{j=1}^N$ and the
precise definitions of solutions to \eqref{equax}.

\begin{enumerate}
\item[\bf (H1).] $e_j\in C_b^{\infty}(\mathbb{R}^d)$ such that
\begin{equation*}
    \lim\limits_{|\xi|\to\infty} \zeta(\xi)|\partial^{\g}e_j(\xi)|=0,
\end{equation*}
where $\g$ is a multi-index such that $|\g|\leq 3$, $1\leq j \leq N$
and
\begin{equation*}
    \zeta(\xi) = \left\{
       \begin{array}{ll}
         1+|\xi|^2, & \hbox{if $d \neq 2$;} \\
         (1+|\xi|^2)(\ln (3+|\xi|^2))^2, & \hbox{if $d=2$.}
       \end{array}
     \right.
\end{equation*}
\end{enumerate}

\begin{definition} \label{defx}
Let $x\in H^1$ and let $\a$ satisfy $1<\a<\9$ if $d=1,2$ or
$1<\a\leq 1+\frac{4}{d-2}$ if $d\geq 3$. Fix $0<T<\9$.

A strong solution of \eqref{equax} is a pair $(X,\tau)$, where
$\tau(\leq T)$ is an $(\mathscr{F}_t)$-stopping time, and
$X=(X(t))_{t\in[0,T]}$ is an $H^1$-valued continuous
$(\mathscr{F}_t)$-adapted process, such that $|X|^{\alpha-1}X\in
L^1(0,\tau;H^{-1})$, $\mathbb{P}-a.s$, and it satisfies
$\mathbb{P}-a.s$
\begin{align} \label{weakx}
    X(t)
    =&x-\int_0^t(i\Delta X(s)+\mu X(s)+\lambda
    i|X(s)|^{\alpha-1}X(s))ds \nonumber \\
     &  +\int_0^t X(s)dW(s),\ \ t\in [0,\tau],
\end{align}
as an equation in $H^{-1}$.

We say that uniqueness holds for \eqref{equax}, if for any two
strong solutions $(X_i, \tau_i)$, $i=1,2$, it holds
$\mathbb{P}$-a.s. that $X_1=X_2$ on $[0,\tau_1\wedge \tau_2]$.
\end{definition}

We refer to \cite{PR07} for the general theory of infinite
dimensional stochastic equations. It is easy to check that,
$\int_0^t
X(s)dW(s)$ in Definition \ref{defx} is an $H^1$-valued continuous stochastic integral.\\

The main global well-posedness result in this paper is as follows.

\begin{theorem}\label{thmx}
 Assume $(H1)$. Let $\a$ satisfy
(\ref{a-defo}) and (\ref{a-fo}) in the defocusing and focusing cases
respectively. Then for each $x\in H^1$ and $0<T<\infty$, there
exists a unique strong solution $(X,T)$ of (\ref{equax}) in the
sense of Definition \ref{defx}, such that
\begin{align} \label{thmx1}
   X\in L^2(\Omega; C([0,T];H^1)) \cap L^{\a+1}(\Omega;
   C([0,T];L^{\a+1})),
\end{align}
and
\begin{align}\label {thmx2}
  X \in L^\g(0,T;W^{1,\rho}),\ \ \mathbb{P}-a.s.,
\end{align}
where $(\rho,\g)$ is any Strichartz pair (see Lemma \ref{Stri-S}
below).

Furthermore, for $\mathbb{P}-a.e$ $\omega$, the map $x \to
X(\cdot,x,\omega)$ is continuous from $H^1$ to $C([0,T];H^1)\cap
L^\g(0,T;W^{1,\rho})$.
\end{theorem}

The key approach here (as in \cite{BRZ14}) is based on the rescaling
transformation that reduces the stochastic equation (\ref{equax}) to
a random Schr\"odinger equation (see (\ref{equay})), to which one
can apply the sharp deterministic estimates, e.g. the Strichartz
estimates established in \cite{MMT08} for large perturbations of the Laplacian. \\

This paper is structured as follows. In Section \ref{LWP} we
establish the local existence, uniqueness and blowup alternative of
solutions to equation (\ref{equax}). Then in Section \ref{Hami-Pri}
we derive a priori estimates of the energy from the Hamiltonian,
which lead to the global well-posedness in the subcritical case in
Section \ref{mainthm}. An important role in our proofs is played by
It\^o's formulae for the $L^p$- and $H^1$- norms, which can be
heuristically computed very easily. The rigorous proofs are much
harder and are contained in Section \ref{Ito-LpH1}. Furthermore,
some technical proofs are postponed to the
Appendix, i.e. Section \ref{App-proof},  for simplicity of exposition. \\

{ \it  \bf Notations.} For $1\leq p \leq \9$,
$L^p=L^p(\mathbb{R}^d)$ is the space of all $p$-integrable complex
valued functions with the norm $|\cdot|_{L^p}$. $L^q(0,T;L^p)$
denotes the measurable functions $u:[0,T]\to L^p$ such that $t\to
|u(t)|_{L^p}$ belongs to $L^q(0,T)$. $C([0,T];L^p)$ similarly
denotes the continuous $L^p$-valued functions with the $\sup$ norm
in $t$.

As usual, $W^{1,p}=W^{1,p}(\mathbb{R}^d)$ is the classical Sobolev
space, i.e. $W^{1,p}=\{u\in L^p: \na u \in L^p\}$ with the norm
$\|u\|_{W^{1,p}}=|u|_{L^p}+|\na u|_{L^p}$. Here $\na
=(\partial_{1},...,\partial_{d})$ with
$\partial_{k}:=\frac{\partial}{\partial{x_k}}$, $1\leq k\leq d$.
Moreover, the spaces $L^q(0,T;W^{1,p})$ and $C([0,T];W^{1,p})$ are
understood similarly as above. We also use the notation
$\partial^{\gamma} =
   \partial_1^{\gamma_1}\cdot\cdot\cdot\partial_d^{\gamma_d}$
for any multi-index $\gamma=(\gamma_1,...,\gamma_d)$ with $\gamma_j
\in\mathbb{N}$. The order of $\gamma$ is $|\gamma| = \gamma_1 +
\cdot\cdot\cdot + \gamma_d $, and if $|\gamma|=0$,
$\partial^{\gamma} f = f$.

In the special case $p=2$, $L^2$ is the Hilbert space endowed with
the scalar product
$$\<u,v\>=\int\limits_{\mathbb{R}^d}u(\xi)\overline{v}(\xi)d\xi;~u,v\in L^2.$$
For simplicity, we set $|\cdot|_{2}=|\cdot|_{L^2}$. Let
$H^1=W^{1,2}$ and $H^{-1}$ be the dual space of $H^1$. Their norms
are denoted by $|\cdot|_{H^k}$, $k=\pm 1$.

$C_c^{\9}(\mathbb{R}^d)$ denotes the compactly supported smooth
functions on $\mathbb{R}^d$. We use $\mathcal{S}$ and $\mathcal{S}'$
for the rapidly decreasing functions and the tempered distributions
respectively. Then for $f\in\mathcal{S}$, $\widehat{f}$ means the
Fourier transform, i.e. $\wh f (\eta) = \int f(\xi) e^{-i\xi\cdot
\eta} d\xi$, and for $f\in \mathcal{S}'$, $f^{\vee}$ denotes the
inverse Fourier transform of $f$, i.e. $f^{\vee}(\xi) =
\frac{1}{(2\pi)^d} \int f(\eta) e^{i\xi\cdot \eta} d\eta$.

We use $C$, $\widetilde{C}$ for various constants that may change
from line to line.

\section{Local results} \label{LWP}

In this section, we will establish the local existence, uniqueness
and blowup alternative for equation (\ref{equax}). The main result
is given in Theorem \ref{localx} below.

\begin{theorem} \label{localx}
Assume $(H1)$. Let $\a$ satisfy $1<\a<\9$ if $d=1,2$, or, $1<\a\leq
1+\frac{4}{d-2}$ if $d\geq 3$. For each $x\in H^1$ and $0<T<\9$,
there is a sequence of strong solutions $(X_n,\tau_n)$ of
\eqref{equax}, $n\in\mathbb{N}$, where $\tau_n$ is a sequence of
incresing stopping times, and uniqueness holds in the sense of
Definition \ref{defx}. For every $n\geq 1$, it holds
$\mathbb{P}$-a.s that
\begin{equation}
    X_n|_{[0,\tau_n]}\in C([0,\tau_n];H^1)\cap L^\g(0,\tau_n;W^{1,\rho}),
\end{equation}
where $(\rho,\g)$ is any Strichartz pair.

Moreover, defining $\tau^*(x)=\lim\limits_{n\to \9} \tau_n$ and
$X=\lim\limits_{n\to \9} X_n \chi_{[0,\tau^*(x))}$, we have the
blowup alternative, that is, for $\mathbb{P}$-a.e $\omega$, if
$\tau_n(\oo)<\tau^*(x)(\oo)$, $\forall n \in \mathbb{N}$, then
\begin{equation} \label{subcritical-a}
    \lim\limits_{t\to \tau^*(x)(\oo)}|X(t)(\oo)|_{H^1}=\9,\ \ if\
    1<\a<1+\frac{4}{(d-2)_+},\ d\geq 1,
\end{equation}
and
\begin{equation}\label{critical-a}
    \|X(\oo)\|_{L^{\frac{2(d+2)}{d-2}}(0,\tau^*(x)(\oo);L^{\frac{2(d+2)}{d-2}})}=\9,\
    \ if\ \a=1+\frac{4}{d-2},\ d\geq 3.
\end{equation}
\end{theorem}

\begin{remark} \label{localx-remark}
As seen below in the proof of Proposition \ref{localy} if the norm
in \eqref{subcritical-a} or \eqref{critical-a} is finite
$\mathbb{P}$-a.s., then $\tau^*(x)=T$, $\mathbb{P}$-a.s.
\end{remark}

The key tool to prove Theorem \ref{localx} is based on the rescaling
approach as used in \cite{BRZ14}. Namely, we apply the rescaling
transformation
\begin{align} \label{rescal}
    X=e^Wy
\end{align} to reduce the original stochastic equation
(\ref{equax}) to the random Schr\"odinger equation
\begin{equation}\label{equay}
    \barr{l}
     \dd  \frac{\partial y(t,\xi)}{\partial t} = A(t)y(t,\xi) - \lambda i
    e^{(\alpha-1)Re W(t,\xi)}|y(t,\xi)|^{\alpha-1}y(t,\xi), \vsp
    y(0)=x. \earr
\end{equation}
Here
\begin{equation} \label{Op-A}
A(t)y(t,\xi) := -i (\Delta + b(t,\xi)\cdot \nabla + c(t,\xi))
y(t,\xi)
\end{equation}
with $b(t,\xi)=2 \nabla W(t,\xi)$, $c(t,\xi)= \sum\limits_{j=1}^d
(\partial_j W(t,\xi))^2 + \Delta
W(t,\xi)-i(\mu(\xi)+\widetilde{\mu}(\xi))$ and
$\widetilde{\mu}(\xi)=\frac{1}{2}\sum\limits_{j=1}^N \mu_j^2
e_j^2(\xi).$\\

Analogously to Definition \ref{defx}, the solutions to \eqref{equay}
are defined as follows.
\begin{definition} \label{defy}
Let $x\in H^1$, $0<T<\9$, and  $\a \in (1,\9)$ if $d=1,2$, or
$\a\in(1,1+\frac{4}{d-2}]$ for $d\geq 3$.  The strong solution
$(y,\tau)$ and uniqueness of \eqref{equay} are defined similarly as
in Definition \ref{defx}, just with the modifications that $X$ and
\eqref{weakx} are replaced, respectively, by $y$ and the equation
\begin{align} \label{weaky}
   y(t) = x + \int_0^t A(s)y(s)ds -\int_0^t \lbb i e^{(\a-1)ReW(s)}
   |y(s)|^{\a-1} y(s) ds.
\end{align}
\end{definition}

\begin{remark}
The equivalence between two strong solutions $(X,\tau)$ and
$(y,\tau)$ of $(\ref{equax})$ and $(\ref{equay})$, respectively, can
be proved similarly as in the proof of Lemma $6.1$ in \cite{BRZ14}.
We also refer to \cite{Z14} for more details.
\end{remark}

Therefore, it is equivalent to prove the local results for the
random equation \eqref{equay}. We have the following

\begin{proposition} \label{localy}
Assume the conditions in Theorem \ref{localx} to hold. For each
$x\in H^1$ and $0<T<\9$, there is a sequence of strong solutions
$(y_n,\tau_n)$ of \eqref{equay}, $n\in\mathbb{N}$, where $\tau_n$ is
a sequence of incresing stopping times, and uniqueness holds in the
sense of Definition \ref{defy}. For every $n\geq 1$, it holds
$\mathbb{P}$-a.s that
\begin{equation} \label{locthm1}
    y_n|_{[0,\tau_n]}\in C([0,\tau_n];H^1)\cap L^\g(0,\tau_n;W^{1,\rho}),
\end{equation}
where $(\rho,\g)$ is any Strichartz pair.

Moreover, defining $\tau^*(x)=\lim\limits_{n\to \9}\tau_n$ and
$y=\lim\limits_{n\to \9} y_n \chi_{[0,\tau^*(x))}$, we have the
blowup alternative, namely, for $\mathbb{P}$-a.e $\omega$ if
$\tau_n(\oo)<\tau^*(x)(\oo)$, $\forall n \in \mathbb{N}$, then
\begin{equation*}
    \lim\limits_{t\to \tau^*(x)(\oo)}|y(t)(\oo)|_{H^1}=\9,\ \ if\
    1<\a<1+\frac{4}{(d-2)_+},\ d\geq 1,
\end{equation*}
and
\begin{equation*}
    \|y(\oo)\|_{L^{\frac{2(d+2)}{d-2}}(0,\tau^*(x)(\oo);L^{\frac{2(d+2)}{d-2}})}=\9,\
    \ if\ \a=1+\frac{4}{d-2},\ d\geq 3.
\end{equation*}
\end{proposition}

Inspired by the deterministic case, the local well-posedness of
\eqref{equay} depends crucially on the dispersive properties of the
linear part in \eqref{equay}. Hence, in order to prove Proposition
\ref{localy}, let us first introduce the evolution operators and
Strichartz estimates in Sobolev spaces.

\begin{lemma} \label{evolu-op}
For $\mathbb{P}-a.e.\omega$, the operator $A(t)$ defined in
$(\ref{Op-A})$ generates evolution operators $U(t,s)=U(t,s,\omega)$
in the space $H^1(\mathbb{R}^d)$, $0\leq s\leq t\leq T$. Moreover,
for each $x\in H^1(\mathbb{R}^d)$ and $s\in[0,T]$, the process
$[s,T]\ni t \to U(t,s)x$ is continuous and
$(\mathscr{F}_t)$-adapted, hence progressively measurable with
respect to the filtration $(\mathscr{F}_t)_{t\geq s}$.
\end{lemma}

{\it \bf Proof.}  This lemma is based on \cite{D69} and can be
proved analogously as Lemma $3.3$ in \cite{BRZ14} (see also
\cite{Z14}). \hfill $\square$

\begin{lemma} \label{Stri-S}
Assume $(H1)$. Then for any $T>0$, $u_0\in H^1$ and $f\in
L^{q_2'}(0,T;W^{1,p_2'})$, the solution of
\begin{equation} \label{stri2}
    u(t)=U(t,0)u_0 + \int_0^t U(t,s) f(s) ds, 0\leq t \leq T,
\end{equation}
satisfies the estimates
\begin{equation} \label{stri-l}
    \|u\|_{L^{q_1}(0,T;L^{p_1})}\leq
    C_T(|u_0|_{2}+\|f\|_{L^{q_2'}(0,T;L^{p_2'})}),
\end{equation}
and
\begin{equation} \label{stri-s}
    \|u\|_{L^{q_1}(0,T;W^{1,p_1})}\leq
    C_T(|u_0|_{H^1}+\|f\|_{L^{q_2'}(0,T;W^{1,p_2'})}),
\end{equation}
where $(p_1,q_1)$ and $(p_2,q_2)$ are Strichartz pairs, namely
\begin{equation*}
    (p_i,q_i)\in[2,\infty] \times
    [2,\infty]: \frac{2}{q_i}=\frac{d}{2}-\frac{d}{p_i},~if~d
    \neq 2,
\end{equation*}
or
\begin{equation*}
    (p_i,q_i)\in[2,\infty) \times
    (2,\infty]:\frac{2}{q_i}=\frac{d}{2}-\frac{d}{p_i},~if~d
    =2,
\end{equation*}
Furthermore, the process $C_t$, $t\geq 0$, can be taken to be
$(\mathscr{F}_t)$-progressively measurable, increasing and
continuous.
\end{lemma}
(See the Appendix for the proof.) \\

{\it \bf  Proof of Proposition \ref{localy}.} It is equivalent to
solve the weak equation $(\ref{weaky})$ in the mild sense, namely
\begin{equation} \label{mildy}
     y=U(t,0)x-\lambda i\int_0^tU(t,s) e^{(\alpha-1)ReW(s)}
    g(y(s))ds,
\end{equation}
where $g(y)=|y|^{\a-1}y$. The following fixed point arguments are
standard in the deterministic case (see e.g. \cite{K89} and
\cite{LP09}). However, we emphasize that we have to secure the
$(\mathscr{F}_t)$-adaptedness of the solutions, which allows us
later to apply It\^{o}'s formula to obtain a priori estimates (see
also \cite{BRZ14}). \\

Let us first consider the case $d\geq 3$. Choose the  Strichartz
pair
$(p,q)=(\frac{d(\alpha+1)}{d+\alpha-1},\frac{4(\alpha+1)}{(d-2)(\alpha-1)})$,
set $\mathcal {X}=C([0,T];L^2)\cap L^q(0,T;L^p)$, $\mathcal
{Y}=C([0,T];H^1)\cap L^q(0,T;W^{1,p})$, and consider the integral
operator
\begin{equation} \label{op-F}
    F(y)(t)=U(t,0)x-\lambda i\int_0^tU(t,s)( e^{(\alpha-1)ReW(s)}
    g(y(s)))ds,\ t\in[0,T],
\end{equation}
defined for $y\in \mathcal{Y}$.

We claim that
 \begin{align} \label{fyy}
 F(\mathcal{Y})\subseteq \mathcal{Y}.
\end{align}

In fact, by the Strichartz estimates in Lemma \ref{Stri-S}
\begin{equation} \label{fyy.1}
    \|F(y)\|_{L^q(0,T;W^{1,p})}
    \leq
    C_T\left[|x|_{H^1}
    + \|e^{(\alpha-1)ReW}g(y)\|_{L^{q'}(0,T;W^{1,p'})} \right].
\end{equation}
To estimate the right-hand side, we have that
\begin{align} \label{fyy.2}
    & \| e^{(\alpha-1)ReW} g(y)\|_{L^{q'}(0,T;W^{1,p'})} \nonumber
    \\
    \leq& D_1(T)\left(\||y|^{\alpha-1}y\|_{L^{q'}(0,T;L^{p'})}
    +\||y|^{\alpha-1}|\nabla y| \|_{L^{q'}(0,T;L^{p'})}\right),
\end{align}
where in the last inequality we have used  $|\nabla g(y)|\leq \alpha
|y|^{\alpha-1} |\nabla y|$,
$|\nabla(e^{(\alpha-1)ReW}g(y))|\leq|e^{(\alpha-1)W}|\left[(\alpha-1)|\nabla
W||g(y)|+|\nabla g(y)|\right]$ and $D_1(T):=\alpha  (|\nabla
W|_{L^{\infty}(0,T;L^{\infty})}+2)e^{(\alpha-1)|W|_{L^{\infty}(0,T;L^{\infty})}}$.

With our choice of $(p,q)$, it is easy to verify that
$(\frac{1}{p'},\frac{\alpha}{q})=(\alpha-1)(\frac{1}{(\alpha-1)l},\frac{1}{q})+(\frac{1}{p},\frac{1}{q})$,
where $\frac{1}{l}=\frac{1}{p'}-\frac{1}{p}$, satisfying
$\frac{1}{(\alpha-1)l}=\frac{1}{p}-\frac{1}{d}$. Hence, from
H\"{o}lder's inequality and the Sobolev imbedding
$|y|_{L^{(\alpha-1)l}} \leq D |y|_{W^{1,p}}$ it follows that
\begin{align} \label{fyy.3}
       \||y|^{\a-1}y\|_{L^{q'}(0,T;L^{p'})}
       \leq& T^{\theta}
       \||y|^{\a-1}y\|_{L^{\frac{q}{\a}}(0,T;L^{p'})} \nonumber \\
       \leq& T^{\theta} \|y\|^{\alpha-1}_{L^{q}(0,T;L^{(\alpha-1)l})}
       \|y\|_{L^{q}(0,T;L^{p})} \nonumber \\
       \leq& D^{\a-1} T^{\theta}  \|y\|^{\alpha-1}_{L^{q}(0,T;W^{1,p})}
       \|y\|_{L^{q}(0,T;L^{p})},
\end{align}
with $\theta=\frac{1}{q'}-\frac{\alpha}{q}\geq 0$, and also
\begin{align} \label{fyy.4}
\||y|^{\a-1}|\nabla y| \|_{L^{q'}(0,T;L^{p'})} \leq D^{\a-1}
T^{\theta} \|y\|^{\alpha-1}_{L^{q}(0,T;W^{1,p})}
       \|\nabla y\|_{L^{q}(0,T;L^{p})}.
\end{align}

Thus, inserting (\ref{fyy.3}), (\ref{fyy.4}) into (\ref{fyy.2}) and
(\ref{fyy.1}) yields that for $y\in \mathcal{Y}$
\begin{equation}\label{esti1}
    \|F(y)\|_{L^{q}(0,T;W^{1,p})}
    \leq C_T \left[|x|_{H_1}
    +D_2(T) T^{\theta} \|y\|^{\alpha}_{L^{q}(0,T;W^{1,p})}
    \right],
\end{equation}
with $D_2(T)= D_1(T) D^{\a-1}$. Similarly,
\begin{equation}\label{esti2}
    \|F(y)\|_{L^{\infty}(0,T;H^1)}
    \leq C_T \left[|x|_{H_1}
     +D_2(T) T^{\theta}\|y\|^{\alpha}_{L^{q}(0,T;W^{1,p})}
     \right].
\end{equation}
Hence \eqref{esti1} and \eqref{esti2} yield (\ref{fyy}), as claimed. \\

We now start to construct the strong solutions of \eqref{equay} by
similar arguments as in \cite{BRZ14}.

{\it \bf Step $1$.} Fix $\omega\in\Omega$ and consider $F$ on the
set
\begin{align*}
    \mathcal{Y}^{\tau_1}_{M_1}
    =&\{y\in C([0,\tau_1];H^1)\cap L^q(0,\tau_1;W^{1,p}); \\
     &\qquad \sup\limits_{0\leq t \leq\tau_1}|y(t)-U(t,0)x|_{H^1}+\|y\|_{L^q(0,\tau_1;W^{1,p})}\leq
     M_1\},
\end{align*}
where $\tau_1=\tau_1(\oo)\in (0,T]$ and $M_1=M_1(\oo)>0$ are random
variables.

For $y\in\mathcal{Y}^{\tau_1}_{M_1}$ by estimates (\ref{esti1}) and
(\ref{esti2})
\begin{align*}
     \|F(y)-U(\cdot,0)x\|_{L^{\infty}(0,\tau_1;H^1)}+\|F(y)\|_{L^{q}(0,\tau_1;W^{1,p})}
   \leq \ve_1(\tau_1)
     +2C_{\tau_1}D_2(\tau_1) \tau_1^{\theta}M_1^{\alpha},
\end{align*}
where $\ve_1(t):=\|U(\cdot,0)x\|_{L^q(0,t;W^{1,p})}$ is
$(\mathscr{F}_t)$-adapted. By Lemma \ref{Stri-S}, $\ve_1(t)=\|
\chi_{(0,t)}(\cdot)U(\cdot,0)x \|_{L^q(0,T;W^{1,p})}\leq C_T
|x|_{H^1}<\9$, and $\chi_{(0,t)}(\cdot)U(\cdot,0)x \to 0$, as $t \to
0^+$. This implies
\begin{align*}
\ve_1(t)\to 0,\ \ as\ t \to 0^+
\end{align*}

In order to obtain
$F(\mathcal{Y}^{\tau_1}_{M_1})\subset\mathcal{Y}^{\tau_1}_{M_1}$, we
shall choose $M_1$ and $\tau_1$ in such a way that
\begin{equation*}
    \ve_1(\tau_1)
     +2C_{\tau_1}D_2(\tau_1) \tau_1^{\theta}M_1^{\alpha}
    \leq M_1.
\end{equation*}
To this end, we define the real-valued continuous,
$(\mathscr{F}_t)$-adapted process
$$Z_t^{(1)}=2^\a C_t D_2(t) \ve_1^{\a-1}(t) t^{\theta},\ \ t\in[0,T],$$
 choose the $(\mathscr{F}_t)$-stopping time
$$\tau_1=\inf \left\{t\in[0,T],Z_t^{(1)}>\frac{1}{2} \right\} \wedge
T$$ and set $M_1=2\ve_1(\tau_1)$. Then it follows that
$Z_{\tau_1}^{(1)}\leq \frac{1}{2}$ and
$F(\mathcal{Y}^{\tau_1}_{M_1})\subset\mathcal{Y}^{\tau_1}_{M_1}$.

Moreover, the estimates as in the proof of (\ref{esti1}) show that
for $y_1,y_2\in\mathcal{Y}^{\tau_1}_{M_1}$
\begin{align} \label{d-esti}
     &\|F(y_1)-F(y_2)\|_{L^{\infty}(0,\tau_1;L^2)}+\|F(y_1)-F(y_2)\|_{L^{q}(0,\tau_1;L^{p})} \nonumber \\
    \leq& 2C_{\tau_1}
           \|\lambda
            e^{(\alpha-1)ReW}  (g(y_1)-g(y_2))\|_{L^{q'}(0,\tau_1;L^{p'})} \nonumber \\
    \leq& C_{\tau_1}
           D_1(\tau_1)\|(|y_1|^{\alpha-1}+|y_2|^{\alpha-1})|y_1-y_2| \|_{L^{q'}(0,\tau_1;L^{p'})} \nonumber \\
    \leq& C_{\tau_1}
            D_1(\tau_1) D^{\alpha-1} \tau^{\theta}_1 \left(\|y_1\|^{\alpha-1}_{L^{q}(0,\tau_1;W^{1,p})}
            +\|y_2\|^{\alpha-1}_{L^{q}(0,\tau_1;W^{1,p})}\right)\|y_1-y_2\|_{L^{q}(0,\tau_1;L^{p})}  \nonumber  \\
    \leq& 2C_{\tau_1}
           D_2(\tau_1) M_1^{\alpha-1}  \tau^{\theta}_1 \|y_1-y_2\|_{L^{q}(0,\tau_1;L^{p})} \nonumber \\
    \leq&
    \frac{1}{2} \|y_1-y_2\|_{L^{q}(0,\tau_1;L^{p})},
\end{align}
which implies that $F$ is a contraction in $C([0,\tau_1];L^2)\cap
L^q(0,\tau_1;L^{p})$.

Since  $\mathcal{Y}^{\tau_1}_{M_1}$ is a complete metric subspace in
$C([0,\tau_1];L^2)\cap L^q(0,\tau_1;L^{p})$, Banach's fixed point
theorem yields a unique $y\in\mathcal{Y}^{\tau_1}_{M_1}$ with
$y=F(y)$ on $[0,\tau_1]$.

Consequently, setting $y_1(t):=y(t \wedge \tau_1)$, $t\in[0,T]$, and
using similar arguments as in the proof of Step $1$ in Lemma $4.2$
in \cite{BRZ14}, we deduce that $(y_1,\tau_1)$ is a strong solution
of \eqref{equay}, such that $y_1(t)=y_1(t\wedge \tau_1)$,
$t\in[0,T]$, and $y_{1}|_{[0,\tau_1]} \in C([0,\tau_1]; H^1) \cap
L^q(0,\tau_1;W^{1,p})$. \\

{\it \bf Step $2$.} Suppose that at the $n^{th}$ step we have a
strong solution $(y_n,\tau_n)$ of \eqref{equay}, such that
$\tau_n\geq \tau_{n-1}$, $y_n(t)=y_n(t \wedge \tau_n)$, $t\in[0,T]$,
and $y_{n}|_{[0,\tau_n]} \in C([0,\tau_n]; H^1) \cap
L^q(0,\tau_n;W^{1,p})$.

Set
\begin{align*}
    \mathcal{Y}^{\sigma_n}_{M_{n+1}}
    =&\{z\in C([0,\sigma_n];H^1)\cap L^q(0,\sigma_n;W^{1,p}); \\
     &\qquad \sup\limits_{0\leq t\leq\sigma_n}|z(t)-U(t+\tau_n,\tau_n)y_n(\tau_n)|_{H^1}+\|z\|_{L^q(0,\sigma_n;W^{1,p})}\leq
     M_{n+1}\},
\end{align*}
and define the integral operator $F_n$ on $\mathcal{Y}$ by
\begin{align} \label{op-Fn}
    F_n(z)(t)=U(\tau_n+t,\tau_n)y_n(\tau_n)
        -\lambda i\int_0^t &U(\tau_n+t,\tau_n+s)
     e^{(\alpha-1)ReW(\tau_n+s)} g(z(s))ds, \nonumber \\
     &\qquad \qquad t\in[0,T],\ \ z\in\mathcal{Y}.
\end{align}

Analogous calculations as in Step $1$ show that for
$z\in\mathcal{Y}^{\sigma_n}_{M_{n+1}}$
\begin{align*}
     &\|F_n(z)-U(\cdot+\tau_n,\tau_n)y_n(\tau_n)\|_{L^{\infty}(0,\sigma_n;H^1)}+\|F_n(z)\|_{L^{q}(0,\sigma_n;W^{1,{p}})}\\
    \leq& \ve_{n+1}(\sigma_n) + 2 C_{\tau_n+\sigma_n}
    D_2(\tau_n+\sigma_n) \sigma_n^{\theta} M_{n+1}^\a,
\end{align*}
and for $z_1,z_2\in\mathcal{Y}^{\sigma_n}_{M_{n+1}}$
\begin{align*}
      &\|F(z_1)-F(z_2)\|_{L^{\infty}(0,\sigma_n;L^2)}+\|F(z_1)-F(z_2)\|_{L^{q}(0,\sigma_n;L^{p})}\\
     \leq& 2C_{\tau_n+\sigma_n}
           D_2(\tau_n+\sigma_n) M_{n+1}^{\alpha-1} \sigma_n^{\theta}
           \|z_1-z_2\|_{L^{q}(0,\sigma_n;L^{p})}.
\end{align*}
where
$\ve_{n+1}(t):=\|U(\tau_n+\cdot,\tau_n)y_n(\tau_n)\|_{L^q(0,t;W^{1,p})}$
is $(\mathscr{F}_{\tau_n+t})$-adapted and
\begin{align*}
    \ve_{n+1}(t) \to 0,\ \ as\ t\to 0.
\end{align*}

Similarly, we define the continuous
$(\mathscr{F}_{\tau_n+t})$-adapted process
$$Z_t^{(n)}:= 2^\a C_{\tau_n+t} D_2(\tau_n+t) \ve_{n+1}^{\a-1}(t) t^{\theta},~t\in[0,T],$$
 set $$\sigma_n=\inf \left\{t\in[0,T-\tau_n]:Z_t^{(n)}>\frac{1}{2}
\right\} \wedge (T-\tau_n)$$ and choose
$M_{n+1}=2\ve_{n+1}(\sigma_{n})$. It follows that
$F_n(\mathcal{Y}^{\sigma_n}_{M_{n+1}})\subset\mathcal{Y}^{\sigma_n}_{M_{n+1}}$
and $F_n$ is a contraction in $C([0,\sigma_n];L^2)\cap
L^q(0,\sigma_n;L^{p})$. Hence, because
$\mathcal{Y}^{\sigma_n}_{M_{n+1}}$ is a complete metric subspace in
$C([0,\sigma_n];L^2)\cap L^q(0,\sigma_n;L^{p})$,  Banach's fixed
point theorem implies that there is a unique
$z_{n+1}\in\mathcal{Y}^{\sigma_n}_{M_{n+1}}$ such that
$z_{n+1}=F_n(z_{n+1})$ on $[0,\sigma_n]$.

Then, set $\tau_{n+1}=\tau_n+\sigma_n$ and define
\begin{equation*}
    y_{n+1}(t)=\left\{
                 \begin{array}{ll}
                   y_n(t), & \hbox{$t\in [0,\tau_n]$;} \\
                   z_{n+1}((t-\tau_n)\wedge\sigma_n), & \hbox{$ t\in (\tau_n, T]$.}
                 \end{array}
               \right.
\end{equation*}
It follows from the definitions of $F$ and $F_n$ that
$y_{n+1}=F(y_{n+1})$ on $[0,\tau_{n+1}]$, implying
 $y_{n+1}$ is a solution to $(\ref{equay})$ on $[0,\tau_{n+1}]$.
Moreover, using similar arguments as in the proof of Step $2$ in
Lemma $4.2$ and of Lemma $6.2$ in \cite{BRZ14}, we deduce that
$\tau_{n+1}$ is an $(\mathscr{F}_t)$-stopping time and $y_{n+1}$ is
adapted to $(\mathscr{F}_t)$ in $H^1$. Hence, $(y_{n+1},\tau_{n+1})$
is a strong solution of \eqref{equay}, such that
$y_{n+1}(t)=y_{n+1}(t\wedge\tau_{n+1})$, $t\in[0,T]$, and
$y_{n+1}|_{[0,\tau_{n+1}]} \in C([0,\tau_{n+1}]; H^1) \cap
L^q(0,\tau_{n+1};W^{1,p})$.\\

{\it \bf Step $3$.}  Starting from Step $1$ and repeating the
procedure in Step $2$, we finally construct a sequence of strong
solutions $(y_n,\tau_n)$, $n\in\mathbb{N}$, where $\tau_n$ are
increasing stopping times and $y_{n+1}=y_n$ on $[0,\tau_n]$. \\

The integrability property $y\in L^\g(0,\tau_n;W^{1,\rho})$ for any
Strichartz pair $(\rho,\g)$ follows easily from Lemma \ref{Stri-S}
and the estimate $(\ref{esti1})$. \\

To prove the uniqueness, for any two strong solutions $(\wt
y_i,\sigma_i)$, $i=1,2$, define
$\varsigma=\sup\{t\in[0,\sigma_1\wedge \sigma_2]: \wt y_1=\wt y_2\
on\ [0,t]\}$. Suppose that $\mathbb{P}(\varsigma <\sigma_1\wedge
\sigma_2)>0$. For $\omega\in\{\varsigma <\sigma_1\wedge \sigma_2\}$,
we have $\wt y_1(\oo)=\wt y_2(\oo)$ on $[0,\varsigma(\oo)]$  by the
continuity in $H^1$, and for $t\in[0,\sigma_1\wedge
\sigma_2(\oo)-\varsigma(\oo))$
\begin{align*}
   &\|\wt y_1(\oo)- \wt
   y_2(\oo)\|_{L^{q}(\varsigma(\oo),\varsigma(\oo)+t;L^{p})}\\
   \leq& 2  C_{\varsigma(\oo)+t} D_2(\varsigma(\oo)+t) \wt M(t) t^{\theta} \|\wt y_1(\oo)- \wt y_2(\oo)\|_{L^{q}(\varsigma(\oo),\varsigma(\oo)+t;L^{p})},
\end{align*}
where $\wt M(t):= \|\wt
y_1(\oo)\|^{\a-1}_{L^{q}(\varsigma(\oo),\varsigma(\oo)+t;W^{1,p})}+\|\wt
y_2(\oo)\|^{\a-1}_{L^{q}(\varsigma(\oo),\varsigma(\oo)+t;W^{1,p})}
\to 0$ as $t\to 0$. Therefore, with $t$ small enough we deduce that
$\wt y_1(\oo)=\wt y_2(\oo)$ on $[\varsigma(\oo),\varsigma(\oo)+t]$,
hence $\wt y_1(\oo)=\wt y_2(\oo)$ on $[0,\varsigma(\oo)+t]$, which
contradicts the definition of $\varsigma$.\\

Now, we are left with proving the blowup alternative. Let us
consider the subcritical and critical cases respectively.

$(i)$. The subcritical case $1<\a<1+\frac{4}{d-2}$, $\theta>0$:
Suppose that $\mathbb{P}(M^*<\9; \tau_n<\tau^*(x),\forall n\in
\mathbb{N})>0$, where $M^*:=
\sup\limits_{t\in[0,\tau^*(x))}|y(t)|_{H^1}$. Define
$$Z_t:=2^\a (M^*)^{\a-1} C^\a_{T+t}
 D_2(T+t) t^\theta,\ t\in[0,T],$$ and
$$\sigma:=\inf\left\{t\in[0,T]:Z_t>\frac14 \right\}\wedge T.$$
For $\omega\in\{M^*<\9; \tau_n<\tau^*(x),\forall n\in \mathbb{N}\}$,
since $\tau_n(\omega)<T$, $\forall n\in \mathbb{N}$, by the
definition of $\sigma_n$ in Step $2$, we have
$$\sigma_n(\omega)=\inf
\left\{t\in[0,T-\tau_n(\omega)]:Z_t^{(n)}(\omega)>\frac{1}{2}
\right\}.$$

Notice that, for every $n\ge1$, $\ve_{n+1}(t)\leq C_{\tau_n+t}M^*$
due to the Strichartz estimate \eqref{stri-s}. Moreover,
$|y(\tau_n(\omega))|_{H^1}\le M^*$, $C_{\tau_n(\omega)+t}\le
C_{T+t}$ and $D_2(\tau_n(\omega)+t) \leq D_2(T+t)$. It follows that
$Z_t(\omega) \geq Z_t^{(n)}(\omega)$, therefore $\sigma_n(\omega) >
\sigma(\omega)> 0$. Hence
$\tau_{n+1}(\omega)=\tau_n(\omega)+\sigma_n(\omega) >
\tau_n(\omega)+\sigma(\omega)$, which implies $\tau_{n+1}(\omega) >
\tau_1(\omega)+n\sigma(\omega)$ for every $n\ge1$, contradicting the
fact that $\tau_n(\omega) \le T$. Therefore, we have shown the
blow-up
alternative in the subcritical case. \\

$(ii)$. The critical case $\a=1+\frac{4}{d-2}$ with $d\geq 3$,
$\theta =0$: We will adapt the arguments from \cite{C03} and
\cite{CW89}. Set $q_1=\frac{2(d+2)}{d-2}$. Besides the Strichartz
pair $(p,q)=(\frac{2d^2}{d^2-2d+4},\frac{2d}{d-2})$, let us choose
another Strichartz pair $(p_2,p_2)=(2+\frac{4}{d},2+\frac{4}{d})$.
Then $\frac{1}{p_2'}= \frac{\a-1}{q_1} + \frac{1}{p_2}$.

Suppose that $\mathbb{P}(\|y\|_{L^{q_1}(0,\tau^*(x);L^{q_1})}<\9;\
\tau_n<\tau^*(x),\forall n\in\mathbb{N})>0$. For
$\oo\in\{\|y\|_{L^{q_1}(0,\tau^*(x);L^{q_1})}<\9;\
\tau_n<\tau^*(x),\forall n\in\mathbb{N}\}$, we have
$\sigma_n(\oo)=\inf\{t\in[0,T-\tau_n(\oo)]; Z^{(n)}_t(\oo)>\frac
12\}$ and $Z^{(n)}_{\sigma_n(\oo)}(\oo)=\frac 12$. For convenience,
we omit the dependence on $\oo$ below.

From the definition of $F_n$ and the construction of $y$, one can
check that for every $n\geq 1$ and $t\in[0,\tau^*(x)-\tau_n)$
\begin{align*}
  y(\tau_n+t) = U(\tau_n+t,\tau_n)y(\tau_n)
  -\lbb i \int_{\tau_n}^{\tau_n+t} U(\tau_n+t,s) e^{(\a-1)ReW(s)}
  g(y(s)) ds.
\end{align*}
Then by Lemma \ref{Stri-S} and H\"oder's inequality, for every
$n\geq 1$, $t\in[\tau^*(x)-\tau_n)$
\begin{align*}
   \|y\|_{L^{p_2}(\tau_n,\tau_n+t;W^{1,p_2})}
   \leq& C_T|y({\tau_n})|_{H^1} + C_T \|e^{(\a-1)ReW(s)}
  g(y(s))\|_{L^{p_2'}(\tau_n,\tau_n+t;W^{1,p_2'})}\\
  \leq& C_T|y({\tau_n})|_{H^1}
  + C_TD_1(T) \|y\|^{\a-1}_{L^{q_1}(\tau_n,\tau^*(x);L^{q_1})}
  \|y\|_{L^{p_2}(\tau_n,\tau_n+t;W^{1,p_2})}.
\end{align*}

Since $\|y\|_{L^{q_1}(0,\tau^*(x);L^{q_1})}<\9$, we have
$\|y\|_{L^{q_1}(\tau_n,\tau^*(x);L^{q_1})} \to 0$ as $n\to \9$.
Hence, choosing $n$ large enough, such that $C_TD_1(T)
\|y\|^{\a-1}_{L^{q_1}(\tau_n,\tau^*(x);L^{q_1})} <\frac 12$, we have
for $t\in[0,\tau^*(x)-\tau_n)$,
$\|y\|_{L^{p_2}(\tau_n,\tau_n+t;W^{1,p_2})}
   \leq 2 C_T |y({\tau_n})|_{H^1}$,
yielding $$\|y\|_{L^{p_2}(0,\tau^*(x);W^{1,p_2})}<\9.$$ Therefore,
\begin{align*}
   \|y\|_{L^{q}(0,\tau^*(x);W^{1,p})}
   \leq& C_T|x|_{H^1} + C_T \|e^{(\a-1)ReW(s)}
  g(y(s))\|_{L^{p_2'}(0,\tau^*(x);W^{1,p_2'})} \\
  \leq& C_T|x|_{H^1}
  + C_TD_1(T)\|y\|^{\a-1}_{L^{q_1}(0,\tau^*(x);L^{q_1})}  \|y\|_{L^{p_2}(0,\tau^*(x);W^{1,p_2})}<\9.
\end{align*}

Now, we note that for every $n\geq 1$ and $t\in[0,\sigma_n]$
\begin{align*}
   \e_{n+1}(t) = \|U(\tau_n+\cdot,\tau_n)y(\tau_n)\|_{L^{q}(0,t;W^{1,p})}
   \leq  \wt M^*_n +C_T D_2(T) (\wt M^*_n)^{\a},
\end{align*}
where
\begin{align*}
   \wt M^*_n(\oo):=
   \|y(\oo)\|_{L^q(\tau_n(\oo),\tau^*(x)(\oo);W^{1,p})} \to 0,\ as\ n\to \9.
\end{align*}
Then we choose $n$ large enough such that
\begin{align*}
\wt Z^{(n)}(\oo):=2^\a C_T(\oo)D_2(T)(\oo) [\wt M^*_n(\oo) +C_T(\oo)
D_2(T)(\oo) (\wt M^*_n)^{\a}(\oo)]^{\a-1} < \frac 16.
\end{align*}
But this implies $\frac 16 > \wt Z^{(n)}(\oo) >  Z_t^{(n)}(\oo)$ for
any $t\in [0,\sigma_n(\oo)]$, in particular, $\frac 16 > \wt
Z^{(n)}(\oo) >  Z_{\sigma_n(\oo)}^{(n)}(\oo)=\frac 12$, yielding a
contradiction. Therefore, we have proved the blowup alternative in
the critical case and completed the proof of Proposition
\ref{localy} for
the case $d\geq 3$. \\

For the case $d=1,2$, we modify the Strichartz pair $(p,q)$ by
$p=\a+1$ and $q=\frac{4(\a+1)}{d(\a-1)}$. Note that
$(\frac{1}{p'},\frac{1}{q})=(\alpha-1)(\frac{1}{p},0)+(\frac{1}{p},\frac{1}{q})$
and $2< p<\9$. H\"{o}lder's inequality and Sobolev's imbedding
$|y|_{L^{p}}\leq D|y|_{H^1}$ give
\begin{align} \label{esti21.1}
    \||y|^{\a-1}y\|_{L^{q'}(0,T;L^{p'})}
    \leq D^{\a-1} T^{\theta} \|y\|^{\a-1}_{L^{\9}(0,T;H^1)}
    \|y\|_{L^{q}(0,T;L^{p})},
\end{align}
where $\theta=1-\frac{2}{q}>0$, and
\begin{equation} \label{esti21.2}
    \||y|^{\a-1} \nabla y\|_{L^{q'}(0,T;L^{p'})}
    \leq  D^{\a-1} T^{\theta} \|y\|^{\a-1}_{L^{\9}(0,T;H^1)}
\|\nabla y\|_{L^{q}(0,T;L^{p})}.
\end{equation}

Hence the estimates $(\ref{esti1})$ and $(\ref{esti2})$ are
accordingly modified by
\begin{equation}\label{esti21}
    \|F(y)\|_{L^{q}(0,T;W^{1,p})}
    \leq C_T \left[|x|_{H_1}
    +D_2(T) T^{\theta} \|y\|^{\a-1}_{L^{\9}(0,T;H^1)}\|y\|_{L^{q}(0,T;W^{1,p})}
\right],
\end{equation}
and
\begin{equation}\label{esti22}
    \|F(y)\|_{L^{\infty}(0,T;H^1)}
    \leq C_T \left[|x|_{H_1}
    +D_2(T) T^{\theta} \|y\|^{\a-1}_{L^{\9}(0,T;H^1)}\|y\|_{L^{q}(0,T;W^{1,p})}
\right].
\end{equation}
Similarly to \eqref{d-esti}, we get
\begin{align} \label{esti21.3}
   &\|F(y_1)-F(y_2)\|_{L^{\infty}(0,T;L^2)} +
   \|F(y_1)-F(y_2)\|_{L^{q}(0,T;L^{p})} \nonumber \\
   \leq& C_T D_2(T) T^{\theta}
    \( \|y_1\|^{\a-1}_{L^{\infty}(0,T;H^1)} +\|y_2\|^{\a-1}_{L^{\infty}(0,T;H^1)}
    \) \|y_1-y_2\|_{L^{q}(0,T;L^{p})}.
\end{align}
Therefore, similar arguments as those after $(\ref{esti1})$ and
$(\ref{esti2})$ yield the asserted results in the case $d=1,2$. This
completes the proof of Proposition \ref{localy}.  \hfill $\square$ \\

From the blowup alternative in Theorem \ref{localx} and Remark
\ref{localx-remark}, we see that global existence follows from an a
priori estimate for the energy, which will be derived from the
Hamiltonian in the next section.

\section{A priori estimate of the energy} \label{Hami-Pri}

Define the Hamiltonian $H$
\begin{equation} \label{Hamil}
    H(u)=\frac{1}{2}\int|\nabla
    u|^2d\xi-\frac{\lambda}{\alpha+1}\int|u|^{\alpha+1}d\xi,\ \
    u\in\ H^1,
\end{equation}
for $1<\a<1+\frac{4}{(d-2)_+}$, $d\geq 1$. Note that $H$ is well
defined by the Sobolev imbedding theorem.

Let $X$, $\tau^*(x)$ be as in Theorem \ref{localx}. The evolution
formula for $H(X)$ is given in Theorem \ref{Hami-Ito} below.

\begin{theorem} \label{Hami-Ito}
Let $\a$ satisfy \eqref{a-defo}. Set $\phi_j=\mu_j e_j$,
$j=1,...,N$. Then $\mathbb{P}$-a.s
\begin{align*}
     &H(X(t))\\
    =&H(x)
       +\int_0^t Re \<-\nabla(\mu X(s)),\nabla X(s) \>_2ds
       + \frac{1}{2} \sum\limits_{j=1}^N \int_0^t |\nabla (X(s)\phi_j)|_2^2ds        \\
    & -\frac{1}{2}\lambda (\alpha-1)\sum\limits_{j=1}^N \int_0^t \int (Re\phi_j)^2  |X(s)|^{\alpha+1} d\xi ds\\
    & +\sum^N_{j=1} \int_0^t Re \<\nabla(\phi_j
    X(s)),\nabla X(s) \>_2d\beta_j(s)\\
    & -\lambda \sum^N_{j=1} \int_0^t \int Re\phi_j |X(s)|^{\alpha+1} d\xi
    d\beta_j(s),\ \ 0\leq t<\tau^*(x).
\end{align*}
\end{theorem}

\begin{remark}
In the deterministic case $\mu_j=0$, $1\leq j\leq N$, the
Hamiltonian is conserved, i.e. $H(X(t))=H(x)$. In the stochastic
conservative case $\mu_j=-i\widetilde{\mu}_j$, $\widetilde{\mu}_j\in
\mathbb{R}$, $1\leq j\leq N$, the above evolution formula for
$H(X(t))$ coincides with $(4.26)$ in \cite{BD03}.
\end{remark}

{\it \bf  Proof.}  This formula follows heuristically by applying
It\^o's formula to the integrands in $H(X(t))$ with the variable
$\xi$ fixed and then integrating over $\mathbb{R}^d$. But the spaces
$L^2$, $L^p$ consist of equivalent classes of functions, the
delicate problem here is to find a suitable version such that for
every $\xi$ fixed, $(X(t,\xi))_{t\in[0,T]}$ is a continuous
semimartingale, which may not exist.
Therefore, we proceed by approximation to give a rigorous proof. \\

We introduce the operators $\Theta_m$, $m\in \mathbb{N}$, used in
\cite{BD03} and defined for any $f\in \mathcal{S}$ by
\begin{equation*}
    \Theta_mf:=\(\theta(\frac{|\cdot|}{m})\)^{\vee}\ast f\ \
    (=m^d\theta^{\vee}(m\cdot)\ast f),
\end{equation*}
where $\theta\in C_c^{\9}$ is real-valued, nonnegative and
$\theta(x)=1$ for $|x|\leq 1$, $\theta(x)=0$ for $|x|>2$.

By Hausdorf-Young's inequality, since $\int \theta^{\vee} d\xi =1$,
we have for any $p\in [1,\9)$
\begin{align} \label{tha-bdd}
   \|\Theta_m\|_{L^{p} \to L^{p} } \leq C,
\end{align}
where $C=C(p)$ is independent of $m$ and
\begin{align} \label{tha-to}
   \Theta_m f \to f\  \  in\ L^p,\ as\ m\to \9.
\end{align}
Moreover, for any $f\in L^{\frac{\a+1}{\a}}$ we have
\begin{align}
   &\Theta_m f \in L^{\a+1}, \label{tha-lp}\\
   &Re \int i f(\xi)\overline{\Theta_mf(\xi)} d\xi=0. \label{tha-l0}
\end{align}
(See the Appendix for the proof.)

Consider the approximating equation
\begin{equation} \label{app-equx}
     \barr{l} idX_m=\Delta X_mdt-i\mu X_mdt+\lambda
     \Theta_m(g(X_m))dt+iX_mdW,t\in(0,T), \vsp
      X_m(0)=x,  \earr
\end{equation}
where $g(X_m)=|X_m|^{\a-1}X_m$. Since the bound in $(\ref{tha-bdd})$
is independent of $m$, the arguments in the proof of Proposition
\ref{localy} show that there exist unique strong solutions
$(X_{m,n},\tau_n)$ of \eqref{app-equx}, $n\in\mathbb{N}$, where
$\tau_n$ are increasing stopping times, independent of $m$. Define
\begin{align} \label{app-equx.sol}
   X_m : =\lim\limits_{n\to \9} X_{m,n}
\chi_{[0,\tau^*(x))}
\end{align}
with $\tau^*(x)=\lim\limits_{n\to \9}\tau_n$. Set
$q=\frac{4(\a+1)}{d(\a-1)}$. We have $\mathbb{P}-a.s.$
\begin{equation} \label{int-xm}
   R(t):=\sup\limits_{m\geq
1}(\|X_m\|_{C([0,t];H^1)}+\|X_m\|_{L^{q}(0,t;W^{1,\a+1})})<\9,\ \
t<\tau^*(x).
\end{equation}

Moreover, it follows from Lemma \ref{Ito-Lp} and Lemma \ref{Ito-H1}
in Section \ref{Ito-LpH1} that
\begin{align} \label{Ito-m}
     &H(X_m(t)) \nonumber \\
    =& H(x) + \int_0^t Re \<-\nabla(\mu X_m),\nabla(X_m) \>_2 dt
       + \frac{1}{2} \sum\limits_{j=1}^N  \int_0^t |\nabla (X_m(s) \phi_j)|_2^2ds    \nonumber     \\
    & -\frac{1}{2}\lambda (\alpha-1) \sum\limits_{j=1}^N \int_0^t  \int (Re\phi_j)^2  |X_m(s)|^{\alpha+1} d\xi ds  \nonumber  \\
    & -\lambda \int_0^t Re \int i\na [(\Theta_m-1)g(X_m)] \na \overline{X_m} d\xi ds   \\
    & +\sum^N_{j=1} \int_0^t Re \<\nabla(\phi_j
    X_m(s)),\nabla(X_m(s)) \>_2d\beta_j(s)   \nonumber \\
    & -\lambda \sum^N_{j=1} \int_0^t \int Re\phi_j |X_m(s)|^{\alpha+1} d\xi d\beta_j(s). \nonumber
\end{align}

In order to pass to the limit in \eqref{Ito-m}, we note that
$\mathbb{P}$-a.s. for $t<\tau^*(x)$
\begin{equation} \label{xm-x}
   X_m\to X,~in~L^{\9}(0,t;H^1) \cap L^{q}(0,t;W^{1,\a+1})
\end{equation}
(see Section \ref{Ito-LpH1} for the proof).

Let us consider the fifth term in the right hand side of
\eqref{Ito-m} for example. We will show that $\mathbb{P}$-a.s. for
$t<\tau^*(x)$
\begin{align} \label{xm-x.1*}
   \lbb \int_0^t Re \int i\na [(\Theta_m-1)g(X_m)] \na \overline{X_m} d\xi ds
   \to 0,~as~m\to \9.
\end{align}
Indeed, because of \eqref{xm-x} it suffices to show that
$\mathbb{P}$-a.s.
\begin{align}
   \na [(\Theta_m-1)g(X_m)] \to 0, \ in\
   L^{q'}(0,t;L^{\frac{\a+1}{\a}})
\end{align}
for $t<\tau^*(x)$. We note that by \eqref{tha-bdd}
\begin{align*}
   &\|\na [(\Theta_m-1)g(X_m)]\|_{L^{q'}(0,t;L^{\frac{\a+1}{\a}})}\\
   \leq& \|(\Theta_m-1)(\na g(X_m)-\na g(X))\|_{L^{q'}(0,t;L^{\frac{\a+1}{\a}})}
    +\|(\Theta_m-1)\na g(X)\|_{L^{q'}(0,t;L^{\frac{\a+1}{\a}})}\\
   \leq& C \|\na g(X_m)-\na g(X)\|_{L^{q'}(0,t;L^{\frac{\a+1}{\a}})}
    +\|(\Theta_m-1)\na g(X)\|_{L^{q'}(0,t;L^{\frac{\a+1}{\a}})},
\end{align*}
where $C$ is independent of $m$. Using the arguments after
\eqref{condep.h1.2} below we deduce that the first term tends to
$0$. Moreover, the second term also converges to $0$, due to
\eqref{tha-to} and \eqref{tha-bdd}. Therefore, we obtain
\eqref{xm-x.1*}, as claimed.

One easily verifies that we can also take the limit for the
remaining terms in (\ref{Ito-m}) using $(\ref{xm-x})$. Consequently,
we complete the proof. \hfill $\square$ \\

We next prove the a priori estimate of the energy in Theorem
\ref{priesti} below. Before that, let us first state and prove some
technical lemmas.

\begin{lemma} \label{l3}
Let $Y\geq 0$ be a real-valued progressively measurable process. We
have
\begin{equation*}
    \mathbb{E} \(\int_0^tY^2(s)ds \)^{\frac{1}{2}}
    \leq \epsilon \mathbb{E}\sup\limits_{s\leq t}Y(s)
          +C_{\epsilon}\int_0^t\mathbb{E}\sup\limits_{r\leq
          s}Y(r)ds.
\end{equation*}
\end{lemma}
{\it \bf Proof.} This lemma follows easily from the fact that $
\int_0^tY(s)^2ds \leq
 \sup\limits_{s\leq t}Y(s) \int_0^tY(s)ds $ and Cauchy's
 inequality. \hfill $\square$

\begin{lemma} \label{G-N}
For $1<\alpha<1+\frac{4}{d}$, $d\geq 1$, we have
\begin{equation} \label{l2}
     |X|^{\alpha+1}_{L^{\alpha+1}}
     \leq C_{\epsilon}|X|_2^{p}+\epsilon|\nabla X|_2^2,
\end{equation}
where $p>2$.
\end{lemma}

{\it \bf Proof.}  From the Gagliardo-Nirenberg inequality it follows
that  $|X|^{\alpha+1}_{L^{\alpha+1}}\leq C|X|_2^{\beta}|\nabla
    X|_2^{\gamma}$,
where $\beta=(1-\theta)(\alpha+1)$ and
$\gamma=\theta(\alpha+1)\in(0,2)$ with
$\theta=\frac{d(\alpha-1)}{2(\alpha+1)}\in(0,1)$. Then, $(\ref{l2})$
follows immediately from Young's inequality $ab\leq
C_{\epsilon}a^\rho+\epsilon b^\delta$,
$\frac{1}{\rho}+\frac{1}{\delta}=1$, by choosing $\gamma \delta=2$.
\hfill $\square$  \\

Unlike in the conservative case, $|X(t)|_2^2$ is no longer
independent of $t$, but a general martingale (see Lemma $4.3$ in
\cite{BRZ14}). After applying Lemma \ref{G-N} to control
$|X(t)|_{L^{\alpha+1}}^{\alpha+1}$, we also need Lemma \ref{l4}
below to bound the $p$-power of $|X(t)|_2$. Its proof is postponed
to the Appendix.

\begin{lemma}  \label{l4}
Let $p\geq 2$. Then there exists $\widetilde{C}(T)<\infty$ such that
\begin{equation*}
    \mathbb{E}\sup\limits_{t\in [0,\tau^*(x))}|X(t)|_2^p\leq \widetilde{C}(T)<\infty.
\end{equation*}
\end{lemma}

With the above preliminaries, we are now ready to prove the main a
priori estimate for the solution $X$ given by Theorem \ref{localx}.

\begin{theorem} \label{priesti}
Under condition \eqref{a-defo} or \eqref{a-fo}, there exists $\wt
C(T)<\9$, such that
\begin{equation} \label{pri.1}
    \mathbb{E}\left[\sup\limits_{t\in[0, \tau^*(x)) }   \(|\nabla
    X(t)|_2^2 + |X(t)|_{L^{\a+1}}^{\a+1} \) \right]
    \leq \widetilde{C}(T)<\9.
\end{equation}
\end{theorem}

{\it \bf Proof.} $(i)$ First assume that $\lbb=1$. From the
definition of $H$ in $(\ref{Hamil})$ and Theorem \ref{Hami-Ito}, it
follows that $\mathbb{P}$-a.s. for every $n\geq 1$ and $t\in [0,T]$
\begin{align} \label{pri.2*}
     &\frac{1}{2}|\nabla X(t\wedge\tau_n)|_2^2 \nonumber \\
    =&H(x)+\frac{1}{\alpha+1} |X(t\wedge \tau_n)|_{L^{\alpha+1}}^{\alpha+1} \nonumber \\
    & +\int_0^{t\wedge \tau_n} \left[Re \<-\nabla (\mu X(s)),\nabla X(s)
    \>_2
       + \frac{1}{2} \sum\limits_{j=1}^N  |\nabla
       (X(s) \phi_j)|_2^2 \right]ds \nonumber \\
    & -\frac{1}{2}(\alpha-1)\sum\limits_{j=1}^N \int_0^{t\wedge\tau_n} \int (Re \phi_j)^2
    |X(s)|^{\alpha+1} d\xi ds   \\
    & +\sum^N\limits_{j=1} \int_0^{t\wedge\tau_n}Re \<\nabla(\phi_jX(s)),\nabla X(s) \>_2d\beta_j(s) \nonumber \\
    & -\sum^N\limits_{j=1} \int_0^{t\wedge\tau_n} \int Re \phi_j |X(s)|^{\alpha+1} d\xi d\beta_j(s) \nonumber \\
    =&H(x)+\frac{1}{\alpha+1}|X(t\wedge\tau_n)|_{L^{\alpha+1}}^{\alpha+1}
       +J_1(t\wedge\tau_n)+J_2(t\wedge\tau_n)+J_3(t\wedge\tau_n)+J_4(t\wedge\tau_n), \nonumber
\end{align}
where $\tau_n$ is as in Theorem \ref{localx} and $\phi_j=\mu_j e_j$,
$1\leq j\leq N$.

To estimate the second term in the right hand side of
\eqref{pri.2*}, we note that, from (\ref{l2}) and Lemma \ref{l4} it
follows that
\begin{align} \label{pri.3}
     \frac{1}{\alpha+1}\mathbb{E}\sup\limits_{s\leq t\wedge\tau_n}
     |X(s)|_{L^{\alpha+1}}^{\alpha+1}
    \leq&\frac{1}{\alpha+1} C_{\epsilon} \mathbb{E}\sup\limits_{s\leq t\wedge\tau_n}|X(s)|_2^{p}
        + \epsilon \frac{1}{\alpha+1}  \mathbb{E}\sup\limits_{s\leq t\wedge\tau_n}|\nabla X(s)|_2^2
       \nonumber  \\
    \leq& \frac{1}{\alpha+1} C_{\epsilon}\widetilde{C}_T
        + \epsilon \frac{1}{\alpha+1}  \mathbb{E}\sup\limits_{s\leq t\wedge\tau_n}|\nabla X(s)|_2^2.
\end{align}

Concerning $J_1(t \wedge \tau_n)$, we note that
\begin{align*}
    J_1(t)
    \leq  C \int_0^t |\nabla X(s)|_2^2 + |X(s)|_2^2 ds,
\end{align*}
where $C$ depends on $|\phi_j|_{\infty}$ and $| \nabla
\phi_j|_{\infty}$, $1\leq j \leq N$. Hence by Lemma \ref{l4}
\begin{align} \label{pri.4}
    \mathbb{E}\sup\limits_{s\leq t\wedge\tau_n}|J_1(s)|
    \leq C \widetilde{C}(T)t
    + C\int_0^t\mathbb{E}\sup\limits_{r\leq s\wedge\tau_n}|\nabla X(r)|_2^2 d s.
\end{align}

Moreover, since
\begin{equation} \label{pri.4*}
    \mathbb{E}\sup\limits_{s\leq t\wedge \tau_n}|J_2(s)|
    \leq  (\alpha-1)|\mu|_{L^{\infty}}\int_0^t\mathbb{E}\sup\limits_{r\leq s\wedge \tau_n}
    |X(r)|_{L^{\alpha+1}}^{\alpha+1}ds,
\end{equation}
using the estimate (\ref{pri.3}) we have that
\begin{align} \label{pri.5}
    \mathbb{E}\sup\limits_{s\leq t\wedge\tau_n}|J_2(s)|
    \leq& (\alpha-1)|\mu|_{L^{\infty}} C_{\epsilon} \widetilde{C}(T)
    t \nonumber \\
    &+ \epsilon (\alpha-1)|\mu|_{L^{\infty}}  \int_0^t\mathbb{E}\sup\limits_{r\leq s\wedge\tau_n}|\nabla X(r)|_2^2ds.
\end{align}

For $J_3$, the Burkholder-Davis-Gundy inequality yields that
\begin{align*}
    \mathbb{E}\sup\limits_{s\leq t\wedge\tau_n}|J_3(s)|
    &\leq C \mathbb{E}\left[\int_0^{t\wedge\tau_n} \sum^N\limits_{j=1} \(Re\<\nabla (\phi_j X(s)),\nabla X(s)\>_2\)^2ds\right]^{\frac{1}{2}}\\
    &\leq C \mathbb{E} \left(\int_0^{t\wedge\tau_n}|X(s)|_2^4 ds \right)^{\frac{1}{2}}
    +C\mathbb{E} \left(\int_0^{t\wedge\tau_n}|\nabla X(s)|_2^4 ds
    \right)^{\frac{1}{2}},
\end{align*}
where $C$ depends on $|\phi_j|_{\infty}$, $|\nabla
\phi_j|_{\infty}$, $1\leq j \leq N$.  It follows from Lemma \ref{l3}
with $Y$ replaced by $|X(s)|_2^2$ and $|\nabla X(s)|_2^2$
respectively and Lemma \ref{l4} that
\begin{align} \label{pri.6}
    \mathbb{E}\sup\limits_{s\leq t\wedge\tau_n}|J_3(s)|
    \leq&\epsilon C\widetilde{C}(T) +CC_{\epsilon} \widetilde{C}(T)  t
          +\epsilon C\mathbb{E}\sup\limits_{s\leq t\wedge\tau_n}|\nabla
          X(s)|_2^2 \nonumber  \\
     &+CC_{\epsilon}\int_0^t \mathbb{E}\sup\limits_{r\leq s\wedge\tau_n}|\nabla X(r)|_2^2ds.
\end{align}

For the remaining term $J_4$, it follows similarly from  the
Burkholder-Davis-Gundy inequality and Lemma \ref{l3} with $Y$
replaced by $|X|_{L^{\a+1}}^{\a+1}$ that
\begin{align} \label{pri.6*}
    \mathbb{E}\sup\limits_{s\leq t\wedge\tau_n}|J_4(s)|
    \leq& C\mathbb{E}\left[\int_0^{t\wedge\tau_n} \sum^N\limits_{j=1} \left(\int  Re \phi_j |X(s)|^{\alpha+1} d\xi \right)^2ds \right]^{\frac{1}{2}} \nonumber \\
    \leq& C \mathbb{E} \left(\int_0^{t\wedge\tau_n} |X(s)|_{L^{\alpha+1}}^{2(\alpha+1)}ds \right)^{\frac{1}{2}}  \\
    \leq& \epsilon C \mathbb{E}\sup\limits_{s\leq t\wedge\tau_n}
    |X(s)|_{L^{\alpha+1}}^{\alpha+1}
     +CC_{\epsilon} \int_0^t\mathbb{E}\sup\limits_{r\leq s\wedge\tau_n}
    |X(r)|_{L^{\alpha+1}}^{\alpha+1}ds. \nonumber
\end{align}
Then \eqref{pri.3} implies that
\begin{align} \label{pri.7}
    \mathbb{E}\sup\limits_{s\leq t\wedge\tau_n}|J_4(s)|
    \leq& C C_{\epsilon}(\epsilon \widetilde{C}(T)+C_{\epsilon} \widetilde{C}(T) t)
          +\epsilon^2 C \mathbb{E}\sup\limits_{s\leq t\wedge\tau_n}|\nabla
          X(s)|_{2}^2 \nonumber \\
    & +\epsilon  C C_{\epsilon} \int_0^t \mathbb{E}\sup\limits_{r\leq s\wedge\tau_n}|\nabla X(r)|_{2}^2ds.
\end{align}

Now, taking (\ref{pri.3})-(\ref{pri.7}) into (\ref{pri.2*}) and
summing up the respective terms, we conclude that
\begin{align*}
     \frac{1}{2}\mathbb{E}\sup\limits_{s\leq t\wedge\tau_n}|\nabla
     X(s)|_2^2
    \leq& C_1(T)+\epsilon C_2(T)\mathbb{E}\sup\limits_{s\leq t\wedge\tau_n}|\nabla
    X(s)|_2^2 \\
      &+C_3(T) \int_0^t \mathbb{E}\sup\limits_{r\leq s\wedge\tau_n}|\nabla
      X(r)|_2^2ds,
\end{align*}
where the constants $C_k(T)$, $1\leq k\leq 3$, depend on $T$,
$H(x)$, $\a$, $|\phi_j|_{\9}$, $|\nabla\phi_j|_{\9}$,  $1\leq j\leq
N$, and $\mathbb{E}\sup\limits_{t\in[0,\tau^*(x))}|X(t)|_2^p$ with
$p\geq 2$. Then, choosing a sufficiently small $\epsilon$ and using
Gronwall's lemma, we obtain
\begin{equation*}
    \mathbb{E}\sup\limits_{t\in[0,\tau_n]}  |\nabla X(t)|_2^2
    \leq \wt C(T)<\9.
\end{equation*}
Finally, taking $n\rightarrow \infty$ and appylying Fatou's lemma,
we obtain
\begin{align*}
     \mathbb{E}\sup\limits_{t\in[0,\tau^*(x))} |\nabla X(t)|_2^2
    \leq \wt C(T)<\9,
\end{align*}
which implies \eqref{pri.1} by \eqref{l2} and Lemma \ref{l4}.\\

$(ii)$ In the defocusing case $\lbb=-1$, the positivity of the
Hamiltonian simplifies many estimates in the previous case $(i)$,
without using Lemma \ref{G-N} and Lemma \ref{l4}, and the condition
on $\alpha$ is less restrictive.

More precisely, taking (\ref{pri.4}), (\ref{pri.4*}), (\ref{pri.6})
and (\ref{pri.6*}) into Theorem \ref{Hami-Ito} and summing up the
respective terms, we derive that
\begin{align*}
     &\frac{1}{2}\mathbb{E}\sup\limits_{s\leq t\wedge \tau_n}|\nabla X(s)|_2^2
       +\frac{1}{\alpha+1}\mathbb{E}\sup\limits_{s\leq t\wedge \tau_n}|X(s)|_{L^{\alpha+1}}^{\alpha+1}\\
    \leq& C_1(T)+ \epsilon C_2(T)\ \mathbb{E}\sup\limits_{s\leq t\wedge\tau_n}(|\nabla
    X(s)|_2^2+|X(s)|_{L^{\a+1}}^{\a+1})\\
      & +C_3(T) \int_0^t \mathbb{E}\sup\limits_{r\leq s\wedge\tau_n}(|\nabla
      X(r)|_2^2+|X(r)|_{L^{\a+1}}^{\a+1})ds,
\end{align*}
where the constants $C_k(T)$, $1\leq k\leq 3$, depend on $T$,
$H(x)$, $\a$, $|\phi_j|_{\9}$, $|\nabla\phi_j|_{\9}$, $1\leq j\leq
N$, and $\mathbb{E}\sup\limits_{t\in[0,\tau^*(x))}|X(t)|_2^2$.

Therefore, similar arguments as at the end of the previous case
yield (\ref{pri.1}). This completes the proof of Theorem
\ref{priesti}.
\hfill $\square$ \\

\section {Proof of Theorem \ref{thmx}.} \label{mainthm}

By Lemma \ref{l4}, Theorem \ref{priesti} and the fact that
$\|e^{-W}\|_{L^{\9}(0,T;W^{1,\9})}<\9$, $\mathbb{P}$-a.s, it follows
that
\begin{equation}\label{h1y}
    \sup\limits_{0\leq t < \tau^*(x)} |y(t)|_{H^1}^2<\infty,a.s.
\end{equation}

Therefore, $\tau^*(x)=T$, $\mathbb{P}$-a.s, due to the blowup
alternative in Proposition \ref{localy} (see also Remark
\ref{localx-remark}). Modifying the definition of $y$ by $y:=
\lim\limits_{n\to \9} y_n$, we deduce that $(y,T)$ is the unique
strong solution of \eqref{equay}. Therefore, letting $X=e^Wy$, we
conclude that $(X,T)$ is the desired unique strong
solution of \eqref{equax}. \\

The integrability \eqref{thmx1} follows from Lemma \ref{l4} and
Theorem \ref{priesti}, and
\eqref{thmx2} follows from \eqref{locthm1}. \\

It remains to prove the continuous dependence on initial data. Again
it is equivalent to prove this for the random equation
\eqref{equay}, and by Lemma \ref{Stri-S} we only need to show it for
the Strichartz pair $(p,q)=(\a+1,\frac{4(\alpha+1)}{d(\alpha-1)})$.

Suppose that $x_m\rightarrow x$ in $H^1$. Let $(y_m,T)$ be the
unique strong solutions of (\ref{equay}) corresponding to the
initial data $x_m$, $m\geq 1$. Since $|x_m|_{H^1} \leq |x|_{H^1}+1$
for $m\geq m_1$ with $m_1$ large enough, we modify $\tau_1(\leq T)$
in the proof of Proposition \ref{localy} by
\begin{align*}
   \tau_1=\inf\{t\in[0,T]: 2^\a (|x|_{H^1}+1)^{\a-1} C_{t}^\a
D_2(t) t^{\theta} > \frac{1}{2} \} \wedge T,
\end{align*}
such that $\tau_1$ is independent for $m\geq m_1$. Hence, the
contraction arguments there and the uniqueness yield that
\begin{equation*}
    \widetilde{R}:=\sup\limits_{m\geq m_1}(\|y_m\|_{L^{\infty}(0,\tau_1;H^1)}+\|y_m\|_{{L^{q}(0,\tau_1;W^{1,p})}})
    <\9,\ \ \mathbb{P}-a.s.
\end{equation*}

Let us first prove the continuous dependence on initial data on the
interval $[0,\tau_1]$. Analogous calculations as in $(\ref{d-esti})$
show that
\begin{equation} \label{condep1}
     \barr{l}
     \|y_m-y \|_{L^{\infty}(0,t;L^2)}+\|y_m-y \|_{L^{q}(0,t;L^{p})}
     \vsp
    \leq 2C_{T}|x_m-x |_2
         + 2C_{T}D_2(T)\widetilde{R}^{\alpha-1} t^{\theta}  \|y_m-y \|_{L^{q}(0,t;L^{p})}, \earr
\end{equation}
where $\theta=1-\frac{2}{q}>0$. Then taking $t$ small and
independent of $m(\geq m_1)$, we have
\begin{align} \label{condep-lp}
\|y_m-y\|_{L^{\infty}(0,t;L^2)}+\|y_m-y\|_{L^{q}(0,t;L^{p})} \to 0,\
as\ m\to\9.
\end{align}
In particular, it follows that
\begin{align} \label{condep-lp*}
   y_m \to y,\ \ in\ measure\ dt\times d\xi,\ \ as\ m\to \9.
\end{align}

Next, to obtain that
\begin{align} \label{condep-h1}
 \|y_m-y \|_{L^{\infty}(0,t;H^1)}+\|y_m-y \|_{L^{q}(0,t;W^{1,p})}
 \to 0,
\end{align}
we use equation $(\ref{evolu-du-mild})$ in the Appendix to derive
that for $m\geq m_1$
\begin{align} \label{condep.h1.equ}
    \nabla (y_m-y )
    =&U(t,0)\nabla (x_m-x )  + \int_0^tU(t,s)\bigg\{ i(D_j\nabla \widetilde{b}^j +\nabla \widetilde{b}^jD_j +\nabla
        \widetilde{c})(y_m-y ) \nonumber \\
    & \qquad   - \lambda i \nabla \left[ e^{(\alpha-1)ReW(s)}
        \(g(y_m(s))-g(y (s))\)\right]\bigg\}ds,
\end{align}
where $g(y)= |y|^{\a-1}y$.

We note that, by Proposition $2.3 (a)$ in \cite{MMT08} and
\eqref{stri-lp-x} in the Appendix, using a similar estimate as in
\eqref{esti21.1}, we obtain
\begin{align} \label{condep.h1.1}
     &\| i(D_j\nabla \widetilde{b}^j +\nabla \widetilde{b}^jD_j +\nabla \widetilde{c}) (y_m-y ) \|_{\wt X'_{[0,t]}} \nonumber \\
    \leq & \kappa_T \|y_m-y\|_{\wt X_{[0,t]}} \nonumber \\
    \leq & \kappa_T C_T |x_m-x|_2
          +\kappa_T C_T \|e^{(\a-1)Re W} (g(y_m)-g(y))\|_{L^{q'}(0,t;L^{p'})} \nonumber \\
    \leq & C(T)|x_m-x |_2 + C(T) t^{\theta} \|y_m-y \|_{L^{q}(0,t;L^{p})},
\end{align}
where $\theta=1-\frac{2}{q}>0$, $\wt X_{[0,t]}$ is the local
smoothing space defined in \cite{MMT08} and $C(T)$ depends on
$\kappa_T$, $C_T$, $\|e^{(\a-1)W}\|_{L^{\9}(0,T;L^{\9})}$ and $\wt
R$.

Then, applying \eqref{stri-lp-x} to \eqref{condep.h1.equ}, we derive
by \eqref{condep.h1.1} and a similar estimate as in \eqref{esti21.1}
that
\begin{align} \label{condep.h1.2}
    & \|\nabla y_m-\nabla y \|_{L^{\infty}(0,t;L^2)}+\|\nabla y_m-\nabla y \|_{L^{q}(0,t;L^{p})} \nonumber\\
    \leq& 2C_T |\na x_m - \na x|_2
      + 2C_T \| i(D_j\nabla \widetilde{b}^j +\nabla \widetilde{b}^jD_j +\nabla \widetilde{c}) (y_m-y ) \|_{\wt
      X'_{[0,t]}} \nonumber \\
    &+ 2C_T \|\lbb i \na [e^{(\a-1)Re W}
      (g(y_m)-g(y))]\|_{L^{q'}(0,t;L^{p'})}\nonumber \\
    \leq& C(T) |  x_m - x |_{H^1}
    +C(T) t^{\theta} \|y_m-y \|_{L^{q}(0,t;L^{p})} \nonumber \\
    &+C(T) \|\na g(y_m) - \na g(y )\|_{L^{q'}(0,t;L^{p'})},
\end{align}
where $C_T$ depends on $\kappa_T$, $C_T$,
$\|e^{(\a-1)W}\|_{L^{\9}(0,T;W^{1,\9})}$ and $\wt R$.

As regards the last term in the right hand side of
\eqref{condep.h1.2}, we note that $
    \nabla g(y)
    = F_1(y)\nabla y + F_2(y)\nabla\overline{y}$,
where $F_1(y)=\frac{\alpha+1}{2}|y|^{\alpha-1}$ and
$F_2(y)=\frac{\alpha-1}{2}|y|^{\alpha-3}y^2$. Then
\begin{align} \label{xmx.h1.2}
     \nabla g(y_m)- \nabla g(y)
    =& F_1(y_m)[\nabla y_m - \nabla y]
       +[F_1(y_m)-F_1(y )] \nabla y  \nonumber \\
    & +F_2(y_m)[\nabla \overline{y_m}- \nabla\overline{y}]
       +[F_2(y_m)-F_2(y)]\nabla \overline{y} \nonumber \\
    =& I_1+I_2+I_3+I_4 .
\end{align}
Since $|I_1|+|I_3|\leq \alpha |y_m|^{\alpha-1}|\nabla y_m- \na y|$,
(\ref{esti21.2}) yields
\begin{align} \label{xmx.h1.3}
   \|I_1+I_3\|_{L^{q'}(0,t;L^{p'})}
   \leq \a D^{\a-1} \wt R^{\a-1}
   t^{\theta}\|y_m-y\|_{L^{q}(0,t;W^{1,p})}.
\end{align}

Thus plugging (\ref{xmx.h1.2}) and (\ref{xmx.h1.3}) into
(\ref{condep.h1.2}), together with (\ref{condep1}), we derive that
\begin{align} \label{condep4}
      &\|y_m-y \|_{L^{\infty}(0,t;H^1)}+\|y_m-y \|_{L^{q}(0,t;W^{1,p})} \nonumber\\
     \leq& C(T)|x_m-x|_{H^1}
           +C(T)  t^{\theta} \|y_m-y \|_{L^{q}(0,t;W^{1,p})} \nonumber\\
      & +C(T) \|I_2+I_4\|_{L^{q'}(0,t;L^{p'})}.
\end{align}

Therefore, choosing $t$ small and  independent of $m(\geq m_1)$, we
deduce that (\ref{condep-h1}) holds once we prove that
\begin{equation}\label{xmx.h1.k24}
    \|I_2+I_4\|_{L^{q'}(0,t;L^{p'})} \to 0,~as~m\to \infty.
\end{equation}

In order to prove $(\ref{xmx.h1.k24})$, by $(\ref{condep-lp})$ we
have for $dt$-a.e. $s\in[0,t]$, as $m\to \9$
\begin{align*}
    |F_1(y_m(s))|_{L^{\frac{p}{p-2}}}
    \to  |F_1(y(s))|_{L^{\frac{p}{p-2}}},
\end{align*}
which, by \eqref{condep-lp*}, implies that for $dt$-a.e. $s\in[0,t]$
\begin{align*}
    F_1(y_m(s))\to F_1(y(s)), \ \ in\
    L^{\frac{p}{p-2}},
\end{align*}
then
\begin{align*}
    [F_1(y_m(s)) - F_1(y(s))]\na y(s) \to 0 , \ \ in\
    L^{p'}.
\end{align*}
Moreover, for  $dt$-a.e. $s\in[0,t]$,
\begin{align*}
   &| [F_1(y_m(s))- F_1(y(s))] \na y(s)|_{L^{p'}}\\
   \leq & \frac{\a+1}{2} D^{\a-1} (\|y_m\|^{\a-1}_{L^{\9}(0,t;H^1)}+\|y\|^{\a-1}_{L^{\9}(0,t;H^1)})
   \|y(s)\|_{W^{1,p}}\\
   \leq & \frac{\a+1}{2} D^{\a-1} \wt R\ \|y(s)\|_{W^{1,p}}\ \   \in
   L^{q'}(0,t),
\end{align*}
Thus, by Lebesgue's dominated convergence theorem, we obtain
\begin{align*}
   \|I_2\|_{L^{q'}(0,t;L^{p'})} \to 0,\ \ as\ m\to \9.
\end{align*}

The proof for $I_4$ is similar. Therefore, we have proved
\eqref{xmx.h1.k24} hence also \eqref{condep-h1} for $t$ small enough
and independent of $m(\geq m_1)$. Reiterating this procedure in
finite steps we obtain (\ref{condep-h1}) on $[0,\tau_1]$.

Now, since $y_m(\tau_1)\to y(\tau_1)$ in $H^1$, similarly we can
extend the above results to $[0,\tau_2]$ with $\tau_2$ depending on
$|y(\tau_1)|_{H^1}$ and $\tau_1\leq \tau_2\leq T$. Reiterating this
procedure, we then obtain increasing stopping times $\tau_n$,
$n\in\mathbb{N}$, depending on $|y(\tau_{n-1})|_{H^1}$, such that
(\ref{condep-h1}) holds on every $[0,\tau_n]$. Finally, as
$\sup\limits_{t\in[0,T]} |y(t)|_2<\9$, $\mathbb{P}$-a.s, using the
proof of the blowup alternative in Proposition \ref{localy}, we
deduce that for $\mathbb{P}$-a.e. $\omega$ there exists
$n(\omega)<\9$ such that $\tau_{n(\omega)}(\omega) = T$. This
implies the continuous dependence on initial data on $[0,T]$ and
consequently completes the proof of Theorem \ref{thmx}. \hfill
$\square$

\section{It\^o-formulae for $L^p$- and $H^1$- norms}
\label{Ito-LpH1}

This section contains the It\^o-formulae for
$|X_m(t)|_{L^{\a+1}}^{\a+1}$ and $|\na X_m(t)|_2^2$, as well as the
asymptotic formula \eqref{xm-x}, which are used in the proof of
Theorem \ref{Hami-Ito} in Section \ref{Hami-Pri}. \\

Let us start with It\^o's formula for $|X_m(t)|_{L^{\a+1}}^{\a+1}$.
First, we note that Theorem $2.1$ in \cite{K10} is not applicable
here, as we do not have $X  \in L^{\a+1}(0,t;W^{1,\a+1}) $ and
$|X|^{\a-1}X \in L^{\a+1}(0,t;L^{\a+1})$ from Theorem \ref{localx}.
However, for the nonlinearity in the approximating equation
\eqref{app-equx}, by \eqref{tha-lp} and \eqref{tha-l0} we have
$\Theta_m(g(X_m)) \in L^{\a+1}$ and $Re \int i g(X_m)
\overline{\Theta_m(g(X_m))} d\xi=0$, which allow to use the
technique from \cite{K10} to obtain the It\^o formula.\\

Let us adapt the same notation from \cite{K10}. Set $h^{\e}=h \ast
\psi_{\e}$ for any locally integrable function $h$ mollified by $
\psi_{\e}$, where $ \psi_{\e}=\e^{-d}\psi(\frac{x}{\e})$ and $\psi
\in C_c^{\9}(\mathbb{R}^d)$ is a real-valued nonnegative function
with unit integral. Recall that $|h^{\e}|_{L^p}\leq |h|_{L^p}$ and
if $h\in L^p$, then $h^{\e}\to h$ in $L^p$ as $\e \to 0$, $p>1$,
which will be used in the later estimates.

\begin{lemma} \label{Ito-Lp}
Let $X_m$ be as in \eqref{app-equx.sol}. Set $p=\a+1$ with
$1<\a<1+\frac{4}{(d-2)_+}$, $d\geq 1$. We have $\mathbb{P}$-a.s.
\begin{align} \label{ito-lp}
  |X_m(t)|_{L^p}^p
  =&|x|_{L^p}^p- p\int_0^t Re \int i\na g(X_m)(s)\na \overline{X_m}(s) d\xi ds \nonumber \\
  &+\frac{1}{2}p(p-2) \sum\limits_{j=1}^N \int_0^t \int
  (Re\phi_j)^2 |X_m(s)|^{p} d\xi ds \\
  &+p \sum\limits_{j=1}^N  \int_0^t  \int Re\phi_j |X_m(s)|^{p} d\xi d\beta_j(s),\ \ 0\leq t<\tau^*(x). \nonumber
\end{align}
Here $g(X_m)=|X_m|^{p-2}X_m$ and $\phi_j=\mu_je_j$, $1\leq j \leq
N$.
\end{lemma}

{\it \bf Proof.} By $(\ref{app-equx})$ we have  $\mathbb{P}$-a.s.
that
\begin{align} \label{ito-lp.0}
  X_m(t) =& x(t) + \int_0^t \left[-i\Delta X_m(s) - \mu X_m(s) - \lbb
  i g_m(s) \right] ds \nonumber \\
   &+ \int_0^t X_m(s) \phi_j d\beta_j(s), \ \ \qquad t<\tau^*(x),
\end{align}
where $g_m(s)=\Theta_m(g(X_m(s)))$, \eqref{ito-lp.0} is considered
as an It\^o equation in $H^{-1}$ and we used the summation
convention over repeated indices for simplicity.

Taking convolution of both sides of \eqref{ito-lp.0} with the
mollifiers $\psi_{\e}$, we have for every $\xi \in \mathbb{R}^d$
that
\begin{align} \label{ito-lp.1*}
    (X_m(t))^{\e}(\xi)&= x^{\e}(\xi)
    + \int_0^t \left[-i\Delta (X_m(s))^{\e}(\xi)- (\mu X_m(s))^{\e}(\xi) - \lbb
  i (g_m(s))^{\e}(\xi) \right] ds  \nonumber \\
  &\qquad + \int_0^t (X_m(s) \phi_j)^{\e}(\xi) d\beta_j(s),\ \
  t<\tau^*(x),
\end{align}
which holds on a set $\Omega_{\xi}\in\mathscr{F}$ with
$\mathbb{P}(\Omega_{\xi})=1$.

In order to find $\wt \Omega\in\mathcal{F}$ with $\mathbb{P}(\wt
\Omega)=1$ such that \eqref{ito-lp.1*} holds on $\wt \Omega$ for all
$\xi\in\mathbb{R}^d$, we need the continuity in $\xi$ of all terms
in  \eqref{ito-lp.1*}. Let us check this for the stochastic integral
term in \eqref{ito-lp.1*}. Set $\sigma_{n,l}=\inf\{s\in[0,\tau_n]:
|X_m(s)|_{H^1}>l\}\wedge \tau_n$. Since the function $\xi \to
(X_m(s)\phi_j)^{\e}(\xi)$ is continuous and
\begin{align*}
 \mathbb{E} \bigg|\sum\limits_{j=1}^N \int_0^{t\wedge\sigma_{n,l}} (X_m(s)\phi_j)^{\e} (\xi)
 d\beta_j(s) \bigg|^2
 \leq (\sum\limits_{j=1}^N |\phi_j|^2_{L^{\9}}) |\psi_{\e}|^2_2\
 l^2 t< \9,
\end{align*}
it follows that $\xi \to \int_0^{t} (X_m(s)\phi_j)^{\e} (\xi)
 d\beta_j(s)$ is continuous on $\{t\leq \sigma_{n,l}\}$. But since
$\sup\limits_{t\in[0,\tau_n]}|X_m(t)|_{H^1}<\9$, $\mathbb{P}$-a.s,
for $\mathbb{P}$-a.e $\omega\in\Omega$ there exists
$l(\omega)\in\mathbb{N}$ such that
$\sigma_{n,l}(\omega)=\tau_n(\omega)$ for all $l\geq l(\omega)$.
Therefore,
\begin{align} \label{K6*.1}
    \bigcup\limits_{l\in\mathbb{N}} \{t\leq \sigma_{n,l}\}
   =\{ t\leq \tau_n\},
\end{align}
implying that $\xi \to \int_0^{t} (X_m(s)\phi_j)^{\e} (\xi)
 d\beta_j(s)$ is continuous on $\{t\leq \tau_n\}$ hence on $\{t <
\tau^*(x)\}$. One can also check the continuity in $\xi$ for the
drift terms in \eqref{ito-lp.1*}.

Therefore, we conclude that \eqref{ito-lp.1*} holds on a full
probability set $\wt \Omega\in\mathscr{F}$ and $\wt \Omega$ is
independent of $\xi\in\mathbb{R}^d$.

Now, we set for simplicity $X^{\e}_m(t)=(X_m(t))^{\e}(\xi)$ and
correspondingly for the respective other terms. Then by It\^{o}'s
formula we have $\mathbb{P}$-a.s.
\begin{align} \label{K6*.1.1}
  |X_m^{\e}(t)|^p=& |x^{\e}|^p
  -p \int_0^t Re(i g(\overline{X_m^{\e}})(s) \Delta
  X_m^{\e}(s))ds
  -p \int_0^t Re(g(\overline{X_m^{\e}})(s) (\mu X_m)^{\e}(s))
  ds \nonumber \\
  &- \lambda p \int_0^t Re (ig(\overline{X_m^{\e}})(s)
  g_m^{\e}(s))ds
  + \frac{p}{2} \int_0^t |X_m^{\e}(s)|^{p-2} |(X_m\phi_j)^{\e}(s)|^2
  ds \nonumber \\
  &+\frac{1}{2}p(p-2)\int_0^t |X_m^{\e}(s)|^{p-4} [Re(\overline{X_m^{\e}}(s)
  (X_m\phi_j)^{\e}(s))]^2 ds \nonumber \\
  &+p\int_0^t Re(g(\overline{X_m^{\e}})(s)
  (X_m\phi_j)^{\e}(s)) d\beta_j(s),\ \ t<\tau^*(x).
\end{align}

We next integrate \eqref{K6*.1.1} over $\mathbb{R}^d$, and it is not
difficult to justify the interchange of integrals by the
deterministic and stochastic Fubini theorem.  We refer to \cite{Z14}
for more details. Therefore, we obtain that
\begin{align} \label{Ito-pp}
  |X_m^{\e}(t)|_{L^p}^p
  =&|x^{\e}|_{L^p}^p
   -p \int_0^t Re \int  i \na g(X_m^{\e})(s) \na \overline{X_m^{\e}}(s) d\xi ds \nonumber \\
  &-p \int_0^t Re \int (\mu X_m)^{\e}(s) g(\overline{X_m^{\e}})(s) d\xi ds \nonumber \\
  &- \lambda p\int_0^t  Re \int ig(\overline{X_m^{\e}})(s) g_m^{\e}(s) d\xi ds  \nonumber  \\
  &+ \frac{p}{2} \int_0^t \int |X_m^{\e}(s)|^{p-2} |(X_m\phi_j)^{\e}(s)|^2 d\xi  ds \nonumber \\
  &+\frac{1}{2}p(p-2)\int_0^t \int |X_m^{\e}(s)|^{p-4}
  [Re(\overline{X_m^{\e}}(s)(X_m\phi_j)^{\e}(s))]^2 d\xi
  ds \nonumber \\
  &+p \int_0^t  Re \int g(\overline{X_m^{\e}})(s)(X_m\phi_j)^{\e}(s)  d\xi d\beta_j(s) \nonumber \\
  =&|x^{\e}|_{L^p}^p+K_1+K_2+K_3+K_4+K_5+K_6.
\end{align}

Now, we can take the limit $\e\to 0$ in $(\ref{Ito-pp})$. Below we
only do that for $K_1$, $K_3$ and $K_6$. The other terms can be
treated similarly.

First, note that as $\e \to 0^+$
\begin{align} \label{Ito-pp.1}
   X_m^{\e} \to X_m, \ \ in\  L^q(0,t;W^{1,p}),
\end{align}
in particularly,
\begin{align} \label{Ito-pp.2}
   X_m^{\e} \to X_m,\ \na X_m^{\e} \to \na X_m  \ \ in\ measure\ dt\times d\xi.
\end{align}

In order to take the limit for $K_1$, it suffices to show that
\begin{align} \label{K1-2}
    \|\na g(X_m^{\e})-\na g(X_m) \|_{L^{q'}(0,t;L^{p'})} \to 0.
\end{align}
To this end, direct calculations show that
\begin{align} \label{K1-3}
    \na g(X_m^{\e})
    = \frac{p-2}{2} |X_m^{\e}|^{p-4} (X_m^{\e})^2 \na \overline{X_m^{\e}}
    + \frac{p}{2} |X_m^{\e}|^{p-2} \na X_m^{\e}.
\end{align}
To treat the first term in the right hand side above, observe that
for $dt-a.e\ s\in[0,t]$ as $\e \to 0$
\begin{align*}
    &||X_m^{\e}|^{p-4}(s) (X_m^{\e})^2(s)|_{L^{\frac{p}{p-2}}}
    =||X_m^{\e}|^{p-2}(s) |_{L^{\frac{p}{p-2}}}
    =|X_m^{\e}(s)|_{L^p}^{p-2}\\
    \to& |X_m(s)|_{L^p}^{p-2}
    =||X_m|^{p-4}(s) (X_m)^2(s)|_{L^{\frac{p}{p-2}}},
\end{align*}
and $
    |\na \overline{X_m^{\e}}(s)|_{L^p} \to  |\na
    \overline{X_m}(s)|_{L^p},$
which yields by \eqref{Ito-pp.2} that, as $\e \to 0$
\begin{equation*}
    \frac{p-2}{2} |X_m^{\e}|^{p-4}(s) (X_m^{\e})^2(s) \na \overline{X_m^{\e}}(s)
    \to \frac{p-2}{2} |X_m|^{p-4}(s) (X_m)^2(s) \na \overline{X_m}(s),
    ~in~L^{p'}.
\end{equation*}
Similar results hold also for the second term in the right hand side
of $(\ref{K1-3})$. Thus for $dt-a.e\ s\in[0,t]$ as $\e \to 0$
\begin{equation} \label{K1-3.1}
   \na g(X_m^{\e})(s) \to \na g(X_m)(s) ,~in ~ L^{p'}.
\end{equation}
Moreover,
\begin{align} \label{K1-3.1.1}
   &|\na g(X_m^{\e})(s) -\na  g(X_m)(s)|_{L^{p'}} \nonumber  \\
   \leq& 2(p-1) |X_m (s) |_{L^p}^{p-2}|\na X_m(s)  |_{L^p}
   \in L^{q'}(0,t),
\end{align}
which implies $(\ref{K1-2})$ by Lebesgue's dominated convergence
theorem. Hence
\begin{align*}
   \lim\limits_{\e\to 0} K_1
   = -p  \int_0^t Re \int i\na g(X_m)(s) \na \overline{X_m}(s) d\xi ds.
\end{align*}

Concerning the term $K_3$ with $g_m^{\e}$ in \eqref{Ito-pp}, first
observe that
\begin{equation*}
   | g(\overline{X_m^{\e}})(s)- g(\overline{X_m})(s)|_{L^{p'}} \to
   0,\ \ |g_m^{\e}(s)-g_m(s)|_{L^p}\to 0,\ \ s\in[0,t],
\end{equation*}
thus as $\e\to 0$
\begin{equation*}
   Re \int ig(\overline{X^{\e}_m})(s)g_m^{\e}(s)  d\xi
   \to  Re \int ig(\overline{X_m})(s)g_m(s)  d\xi.
\end{equation*}
Moreover, by H\"older's inequality, \eqref{tha-lp} and Sobolev's
imbedding theorem we have
\begin{align} \label{K1-3.1.1*}
   \big|Re \int i g(\overline{X_m^{\e}})(s) g_m^{\e}(s) d\xi\big|
   \leq & |g(\overline{X_m^{\e}})(s)|_{L^{p'}} |g_m^{\e}(s)|_{L^p} \nonumber \\
   \leq & C |X_m(s)|_{L^p}^{2(p-1)} \nonumber  \\
   \leq &  C \sup\limits_{s\in[0,t]}|X_m(s)|_{H^1}^{2(p-1)}<\9,
\end{align}
which, by Lebesgue's dominated convergence theorem and
\eqref{tha-l0}, implies that
\begin{align*}
   \lim\limits_{\e\to 0}K_3
   = - \lambda p \int_0^t Re \int i g(\overline{X_m})(s) g_m(s) d\xi
   ds
   =0.
\end{align*}

Finally, as regards the last stochastic term $K_6$ in
\eqref{Ito-pp}, we first prove that for $\sigma_{n,l}$ defined
above, as $\e \to 0$
\begin{equation} \label{k4}
   \mathbb{E}\int_0^{t\wedge\sigma_{n,l}} Re  \left[\int
   g(\overline{X_m^{\e}})(s) (X_m\phi_j)^{\e}(s) d\xi
   - \int g(\overline{X_m})(s) (X_m\phi_j)(s)  d\xi \right]^2 ds \to
   0.
\end{equation}
In fact, using similar arguments as above, we have for
$s\in[0,t\wedge\sigma_{n,l}]$
\begin{equation}
   Re \int g(\overline{X_m^{\e}})(s) (X_m\phi_j)^{\e}(s) d\xi
   -Re \int g(\overline{X_m})(s) (X_m\phi_j)(s)d\xi
   \to 0.
\end{equation}
Furthermore, as in estimate (\ref{K1-3.1.1*}), for
$s\in[0,t\wedge\sigma_{n,l}]$
\begin{equation}
   \big|\int g(\overline{X_m^{\e}})(s) (X_m\phi_j)^{\e}(s)
   d\xi\big|^2
   \leq C \sup\limits_{s\in[0,t\wedge\sigma_{n,l}]}|X_m(s)|^{2p}_{H^1}<Cl^{2p},
\end{equation}
which yields $(\ref{k4})$ by Lebesgue's dominated convergence
theorem. Hence
\begin{align} \label{K6*}
   K_6 \to p \int_0^t \int Re\phi_j |X_m(s)|^p d\xi d\beta_j(s)
\end{align}
in $\mathbb{P}$-measure on $\{t\leq \sigma_{n,l}\}$ as $\e\to 0$,
which implies by \eqref{K6*.1} that \eqref{K6*} holds on $\{ t\leq
\tau_n\}$. Therefore, as $\tau_n\to \tau^*(x)$ $\mathbb{P}$-a.s, we
conclude that \eqref{K6*} holds $\mathbb{P}$-a.s. for $t<\tau^*(x)$.

Therefore, we can pass to the limit $\e \to 0$ in $(\ref{Ito-pp})$.
As $K_2$ and $K_4$ are canceled after
taking the limit, we finally obtain the desired formula $(\ref{ito-lp})$ . \hfill $\square$ \\

Next, we prove the It\^o formula for $|\na X_m|_2^2$.
\begin{lemma} \label{Ito-H1}
Assume the conditions in Lemma \ref{Ito-Lp} to hold. We have
$\mathbb{P}$-a.s. for $ t<\tau^*(x)$
\begin{align} \label{ito-h1}
  |\na X_m(t)|^2_2
  =&|\na x|^2_2
  +2\int_0^t  Re \< -\na (\mu X_m)(s), \na X_m(s)\>_2ds \nonumber\\
  &+ \sum\limits_{j=1}^N\int_0^t |\na(X_m(s) \phi_j)|_2^2 ds
  -2 \lambda\int_0^t  Re \int i \na g_m(s)\na \overline{X_m}(s) d\xi ds\nonumber \\
  &+2\sum\limits_{j=1}^N\int_0^t  Re \< \na (\phi_jX_m(s)),\na X_m(s)\>_2
  d\beta_j(s).
\end{align}
\end{lemma}

{\it \bf Proof.} We follow the ideas from the proof of $(4.14)$ in
\cite{BRZ14} to derive $(\ref{ito-h1})$. Let $\{f_k|k\in
\mathbb{N}\} \subset H^2$ be an orthonormal basis of $L^2$, set
$J_{\e}=(I- \e \Delta)^{-1}$ and $h_{\e}:=J_{\e}(h) \in H^1$ for any
$h\in H^{-1}$. Then we have from equation $(\ref{app-equx})$ that
$\mathbb{P}$-a.s. for $t\in(0,\tau^*(x))$
\begin{equation} \label{equ-xme}
     \barr{l} idX_{m,\e}=\Delta X_{m,\e}dt-i(\mu X_{m})_{\e}dt+\lambda
     g_{m,\e}dt+i(X_{m}\phi_j)_{\e}d\beta_j, \vsp
      X_{m,\e}(0)=x_{\e}, \earr
\end{equation}
where $g_{m,\e}=[\Theta_m(g(X_m))]_{\e}$ and we used the summation
convention.

Since $\partial_l f_k\in H^1$ for each $f_k$, $1\leq l\leq d$,
$k\in\mathbb{N}$, it follows from (\ref{equ-xme}) and Fubini's
theorem that $\mathbb{P}$-a.s. for $t\in(0,\tau^*(x))$
\begin{align*}
   &\<X_{m,\e}(t),\partial_l f_k \>_2 \\
   =&  \<x_{\e}, \partial_l f_k \>_2
   + \int_0^t \<-i \Delta X_{m,\e}(s),\partial_l f_k  \>_2ds
   + \int_0^t \<-(\mu X_{m})_{\e}(s), \partial_l f_k  \>_2 ds\\
   &+ \int_0^t \<-\lambda ig_{m,\e}(s),\partial_l f_k  \>_2 ds
   + \int_0^t \<(X_{m}(s)\phi_j)_{\e}, \partial_lf_k  \>_2 d\beta_j.
\end{align*}
Applying It\^o's product rule and integrating by parts, we deduce
that
\begin{align*}
   &|\<X_{m,\e}(t),  \partial_l f_k \>_2|^2\\
   =&|\<\partial_l x_{\e},   f_k \>_2|^2
   + 2 Re \int_0^t \overline{\<\partial_l X_{m,\e}(s),f_k \>_2}\<-i\partial_l\Delta X_{m,\e}(s), f_k
   \>_2ds \\
   &+2 Re \int_0^t \overline{\< \partial_l X_{m,\e}(s),f_k \>_2}\<-\partial_l (\mu X_{m})_{\e}(s), f_k  \>_2
   ds \\
   &+2 Re \int_0^t \overline{\< \partial_l X_{m,\e}(s),f_k \>_2}\<-\lambda i\partial_l g_{m,\e}(s), f_k  \>_2
   ds \\
   &+2 Re \int_0^t \overline{\< \partial_l X_{m,\e}(s), f_k \>_2}\<\partial_l (X_m(s)\phi_j)_{\e}, f_k  \>_2
   d\beta_j(s)\\
   &+\int_0^t  |\< \partial_l (X_{m}(s)\phi_j)_{\e}, f_k  \>_2|^2
   ds,\ \ t<\tau^*(x),\ \ \mathbb{P}-a.s.
\end{align*}
We note that $\Delta X_{m,\e}$ and $g_{m,\e}$ are in $H^1$, thus the
above integrals make sense. This is the reason why we have
introduced the operator $J_{\e}$.

Now summing over $k\in \mathbb{N}$ and interchanging infinite sum
and integrals (which can be justified easily), we obtain
$\mathbb{P}$-a.s. for all $t\in (0,\tau^*(x))$
\begin{align*}
   & |\partial_l X_{m,\e}(t) |_2^2\\
   =& \sum\limits_{k=1}^{\9} |\<X_{m,\e}(t), \partial_l f_k \>_2|^2\\
   =&|\partial_l x_{\e}|^2_2
    +2 \int_0^t Re \<i\Delta X_{m,\e}(s), \partial_l^2
    X_{m,\e}(s)\>_2 ds \\
    &+2\int_0^t  Re \<-\partial_l (\mu X_m)_{\e}(s),\partial_l X_{m,\e}(s)\>_2ds +\int_0^t |\partial_l(X_m(s) \phi_j)_{\e}|_2^2 ds \nonumber
  \\
  &-2 \lambda\int_0^t  Re \< i \partial_l g_{m,\e}(s), \partial_lX_{m,\e}(s)\>_2 ds
  +2\int_0^t  Re \< \partial_l (X_m(s)\phi_j)_{\e},\partial_l X_{m,\e}(s)\>_2 d\beta_j(s).\nonumber
\end{align*}

Finally, summing over $l:1\leq l\leq d$ and using the fact that
$f_{\e}\to f$ in $H^k$ and $|f_{\e}|_{H^k}\leq |f|_{H^k}$ for
$k=-1,0,1$, we can pass to the limit $\e \to 0$ in the above
equality and then obtain the evolution
formula $(\ref{ito-h1})$. \hfill $\square$\\

We conclude this section with the proof of the asymptotic formula \eqref{xm-x}.\\

{\it \bf Proof of $(\ref{xm-x})$.} This proof is analogous to that
of continuous dependence on initial data in Theorem \ref{thmx},
hence we only give a sketch of it.  Set $q=\frac{4(\a+1)}{d(\a-1)}$.
By the rescaling transformation $X_m=e^Wy_m$, it suffices to prove
that $\mathbb{P}$-a.s. $$y_m\to y,\ \ in~L^{\9}(0,t;H^1) \cap
L^q(0,t; W^{1,\a+1}) ,\  \ t<\tau^*(x).$$

Notice that, $(\ref{app-equx})$ implies that
\begin{align} \label{app-equy-mild}
      y_m=U(t,0)x-\lambda i\int_0^t
           U(t,s)e^{-W(s)}\Theta_m(g(e^{W(s)}y_m(s)))ds.
\end{align}

By \eqref{int-xm}, \eqref{locthm1} and since
$\|W\|_{L^{\9}(0,T;W^{1,\9})}<\9$, we have $\mathbb{P}$-a.s. for $
t<\tau^*(x)$
\begin{align} \label{app-equy*}
   \wt R(t):=& \sup\limits_{m\geq1}
   (\|y_m\|_{C([0,t];H^1)}+\|y_m\|_{L^q(0,t;W^{1,\a+1})}) \nonumber  \\
   &+(\|y\|_{C([0,t];H^1)}+\|y\|_{L^q(0,t;W^{1,\a+1})})<\9.
\end{align}
Moreover, combining $(\ref{mildy})$ and \eqref{app-equy-mild}, we
have
\begin{align} \label{xmx.equ}
y_m-y=-\lambda i\int_0^t
U(t,s)e^{-W(s)}\left[\Theta_m(g(e^{W(s)}y_m(s)))-g(e^{W(s)}y(s))\right]ds.
\end{align}

Now, we first claim that there exists  $t$ small enough and
independent of $m$, such that
\begin{align} \label{xmx.lp}
     \|y_m-y\|_{L^{\9}(0,t;L^2)}+\|y_m-y\|_{L^{q}(0,t;L^{\a+1})}
   \to 0,~as~m\to \9,
\end{align}
in particularly,
\begin{align} \label{xmx.lp*}
   y_m\to y\ \ in\ measure\ dt\times d\xi.
\end{align}

Indeed, applying Strichartz estimate (\ref{stri-l}) to
(\ref{xmx.equ}) we have
\begin{align} \label{xmx.lp.1}
   &\|y_m-y\|_{L^{\9}(0,t;L^2)}+\|y_m-y\|_{L^{q}(0,t;L^{\a+1})} \nonumber \\
   \leq& 2C_T \|e^{-W}\|_{L^{\9}(0,T;L^{\9})}
   \|\Theta_m(g(e^Wy_m))-g(e^Wy)\|_{L^{q'}(0,t;L^{\frac{\a+1}{\a}})} \nonumber \\
   \leq& C(T)
   \|\Theta_m[g(e^Wy_m)-g(e^Wy)]\|_{L^{q'}(0,t;L^{\frac{\a+1}{\a}})}
   + C(T)
   \|(\Theta_m-1)g(e^Wy)\|_{L^{q'}(0,t;L^{\frac{\a+1}{\a}})}  \nonumber \\
   \leq&C(T)t^{\theta}\|y_m-y\|_{L^{q}(0,t;L^{\a+1})}
     + C(T)
   \|(\Theta_m-1)g(e^Wy)\|_{L^{q'}(0,t;L^{\frac{\a+1}{\a}})},
\end{align}
where we used $(\ref{tha-bdd})$ and estimates as in (\ref{esti21.1})
in the last inequality. Here $\theta=1-\frac{2}{q}>0$, $C(T)$
depends on $C_T$, $\|W\|_{L^{\9}(0,T;L^{\9})}$ and $\wt R(t^*)$ with
any fixed $t^*\in(t,\tau^*(x))$. Choosing $t$ small enough and then
using (\ref{tha-to}), we consequently
obtain (\ref{xmx.lp}), as claimed.\\

Next, we prove that for $t$ sufficiently small and independent of
$m$
\begin{align} \label{xmx.h1}
     \|y_m-y\|_{L^{\9}(0,t;H^1)}+\|y_m-y\|_{L^{q}(0,t;W^{1,\a+1})}
   \to 0,~as~m\to \9.
\end{align}
Indeed, from $(\ref{evolu-du-mild})$ in the Appendix it follows that
\begin{align} \label{xmx.h1.equ}
    \nabla (y_m-y)
    =& \int_0^tU(t,s)\bigg\{ i(D_j\nabla \widetilde{b}^j +\nabla \widetilde{b}^jD_j +\nabla
        \widetilde{c})(y_m-y) \nonumber \\
    & \qquad   - \lambda i \nabla \left[ e^{-W}\( \Theta_m(g(e^Wy_m)) - g(e^Wy) \)\right]\bigg\}ds.
\end{align}
Using estimate as in \eqref{condep.h1.1}, together with
\eqref{xmx.lp.1}, we have that
\begin{align} \label{xmx.h1.1}
     &\| i(D_j\nabla \widetilde{b}^j +\nabla \widetilde{b}^jD_j +\nabla \widetilde{c}) (y_m-y) \|_{\wt X'_{[0,t]}} \nonumber \\
    \leq & C(T) t^{\theta} \|y_m-y \|_{L^{q}(0,t;L^{\a+1})}
           + C(T) \|(\Theta_m-1)g(e^Wy)\|_{L^{q'}(0,t;L^{\frac{\a+1}{\a}})},
\end{align}
where $\theta=1-\frac{2}{q}>0$, and  $C(T)$ depends on $\kappa_T$,
$C_T$, $\|W\|_{L^{\9}(0,T;L^{\9})}$ and $\wt R(t^*)<\9$,
$\mathbb{P}$-a.s.

Then, similarly to \eqref{condep.h1.2}, we have for $m\geq 1$
\begin{align} \label{xmx.h1.1*}
    & \|\nabla y_m-\nabla y \|_{L^{\infty}(0,t;L^2)}+\|\nabla y_m-\nabla y \|_{L^{q}(0,t;L^{\a+1})} \nonumber\\
   \leq& C(T)  t^{\theta} \|y_m-y \|_{L^{q}(0,t;L^{\a+1})}
     + C(T) \|\na g(e^Wy_m)-\na g(e^Wy)\|_{L^{q'}(0,t;L^{{\a+1\over
     \a}})} \nonumber \\
    &+C(T)\|(\Theta_m-1)g(e^Wy)\|_{L^{q'}(0,t;W^{1,\frac{\a+1}{\a}})},
\end{align}
where $C(T)$ is independent of $t$ and $m$.

Therefore, applying analogous arguments as those after
\eqref{condep.h1.2} to control the second term, and then using
\eqref{tha-to} to take the limit in the last term, we deduce that
\eqref{xmx.h1} holds for $t$ small enough and independent of $m$.
Reiterating this procedure with estimates as above we conclude
(\ref{xm-x}) for any $t< \tau^*(x)$. \hfill $\square$

\section{Appendix} \label{App-proof}

{\it \bf Proof of Lemma \ref{Stri-S}.} Estimate $(\ref{stri-l})$ is
already proved in Lemma $4.1$ in \cite{BRZ14}. We can use the same
arguments there to derive that
\begin{equation} \label{stri-lp-x}
    \|u\|_{L^{q_1}(0,T;L^{p_1}) \cap \wt X_{[0,T]}}\leq
    C_T(|u_0|_{2}+\|f\|_{L^{q_2'}(0,T;L^{p_2'}) + \wt X'_{[0,T]}}),
\end{equation}
where $\wt X_{[0,T]}$ is the local smoothing space introduced in
\cite{MMT08} up to time $T$ and
$(q_i,p_i)$, $i=1,2$, are Strichartz pairs. \\

Next, we prove the estimate $(\ref{stri-s})$. Since the proof relies
on Theorem $1.13$ and Proposition $2.3$ $(a)$ in \cite{MMT08}, we
adapt the notations there $D_t:=-i
\partial_t$, $D_j:=-i\partial_j$, $1\leq j \leq d$, to rewrite
$(\ref{stri2})$ in the form
$$ D_t u = (D_ja^{jk}D_k+D_j\widetilde{b}^j+\widetilde{b}^jD_j+\widetilde{c})u - if$$ with $a^{jk}=\delta_{jk}$, $\widetilde{b}^j=-i\partial_j
W_t$ and $\widetilde{c}=-\sum\limits_{j=1}^d(\partial_j W)^2
+(\mu+\widetilde{\mu})i$, $1 \leq j,k \leq d$.

Direct computations show
\begin{align} \label{evolu-du}
     D_t \nabla u
    =& (-\Delta+D_j\widetilde{b}^j+\widetilde{b}^jD_j+\widetilde{c}) \nabla u  \nonumber \\
     &+(D_j \nabla\widetilde{ b}^j+\nabla \widetilde{b}^jD_j+ \nabla
    \widetilde{c})u- i \nabla f.
\end{align}
We regard $(\ref{evolu-du-mild})$ as the equation for the unknown
$\nabla u$ and treat the lower order term $(D_j \nabla
\widetilde{b}^j+\nabla \widetilde{b}^jD_j+ \nabla \widetilde{c})u$
as equal terms with  $\nabla f$. This leads to
\begin{align} \label{evolu-du-mild}
   \na u(t) = U(t,0) \na u_0
      + \int_0^t U(t,s)
    \bigg[ i(D_j \nabla\widetilde{ b}^j(s)+\nabla \widetilde{b}^j(s)D_j  + \nabla \widetilde{c}(s))u(s) +\nabla f(s) \bigg] ds.
\end{align}
Hence applying (\ref{stri-lp-x}) to \eqref{evolu-du-mild} and then
using Proposition $2.3$ $(a)$ in \cite{MMT08} to control the lower
order term, we derive that
\begin{align} \label{low-esti}
     &\|\nabla u\| _{L^{q_1}(0,T;L^{p_1}) \cap \wt X_{[0,T]}} \nonumber \\
    \leq& C_T \left[  |\nabla u_0|_2 + \| i( D_j \nabla \widetilde{b}^j+\nabla \widetilde{b}^jD_j + \nabla
    \widetilde{c})u \|_{\wt X'_{[0,T]}}+ \|\nabla f \|_{L^{q'_2}(0,T;L^{p'_2})}  \right]  \nonumber\\
    \leq& C_T \left[  |\nabla u_0|_2 + \kappa_T\|u\|_{\wt X_{[0,T]}} + \|\nabla
    f\|_{L^{q'_2}(0,T;L^{p'_2})}  \right] \\
    \leq& C_T \left[  |\nabla u_0|_2 + C_T\kappa_T(|u_0|_2+\|f\|_{L^{q'_2}(0,T;L^{p'_2})}) + \|\nabla
    f\|_{L^{q'_2}(0,T;L^{p'_2})}  \right] \nonumber \\
    =& C_T(C_T\kappa_T+1) \left[  |u_0|_{H^1} + \|f\|_{L^{q'_2}(0,T;W^{1,p'_2})}  \right], \nonumber
\end{align}
where we also used $(\ref{stri-l})$ to estimate $\|u\|_{\wt
X_{[0,T]}}$ in the last two inequalities. This together with
$(\ref{stri-l})$ yields the
estimate $(\ref{stri-s})$. \\

Now, set
\begin{align}
   C_t=& \sup\{\|U(\cdot,0)u_0\|_{L^{q_1}(0,t;W^{1,p_1})}
   ;|u_0|_{H^1}\le1\} \nonumber \\
   &+\sup\left\{\left\|
   \int^{\cdot}_0U(\cdot,s)f(s)ds\right\|_{L^{q_1}(0,t;W^{1,p_1})};
   \|f\|_{L^{q'_2}(0,t;W^{1,p'_2})}=1\right\}.
\end{align}
Then the asserted properties of $C_t$, $t\geq 0$, follow analogously
as in the proof of Lemma $4.1$ in \cite{BRZ14} (see also
\cite{Z14}). This
completes the proof of Lemma \ref{Stri-S}.     \hfill $\square$\\

{\it \bf Proof of $(\ref{tha-lp})$.} Hausdorf-Young's inequality
shows that
\begin{align*}
     |\Theta_m f|_{L^{\a+1}}=&|(\theta(\frac{|\cdot|}{m}))^{\vee} \ast
     f|_{L^{\a+1}}
     \leq  |(\theta(\frac{|\cdot|}{m}))^{\vee}|_{L^{\frac{\a+1}{2}}}
     |f|_{L^{\frac{\a+1}{\a}}}.
\end{align*}
As $\theta(\frac{|\cdot|}{m}) \in C_c^{\9} \subset \mathcal {S}$,
$(\theta(\frac{|\cdot|}{m}))^{\vee} \in \mathcal {S} \subset
L^{\frac{\a+1}{2}}$, which implies
$|\Theta_m f|_{L^{\a+1}}<\9.$  \hfill $\square$ \\

{\it \bf Proof of $(\ref{tha-l0})$.} For $f\in L^{\frac{\a+1}{\a}}
\cap L^1$, (\ref{tha-l0}) follows from Fourier's inversion formula
and Fubini's theorem. The general case $f\in L^{{\a+1\over \a} }$
follows from a standard approximating procedure.  \hfill $\square$ \\

{\it \bf Proof of Lemma \ref{l4}.} As in the proof of Lemma $4.3$ in
\cite{BRZ14}, we have
\begin{align} \label{l4.1*}
    |X(t)|_2^2=|x|_2^2 +
2\sum\limits_{j=1}^N \int_0^t Re\mu_j <X(s),X(s)e_j>_2 d\beta_j(s),\
\ t<\tau^*(x),\ \mathbb{P}-a.s.,
\end{align}
where $\tau^*(x)$ is as in Theorem \ref{localx}. Then, It\^o's
formula implies
\begin{align*}
    |X(t)|_2^p
    =& |x|_2^p + p\int_0^t|X(s)|_2^{p-2}\sum\limits_{j=1}^N
    Re\mu_j\<X(s),X(s)e_j\>_2d\beta_j(s)\\
    &+ \frac{1}{2} p (p-2)\int_0^t |X(s)|_2^{p-4} \sum\limits_{j=1}^N (Re\mu_j)^2\<X(s),X(s)e_j\>_2^2
    ds, \ \ t<\tau^*(x).
\end{align*}

Hence, by the Burkholder-Davis-Gundy inequality and Lemma \ref{l3}
with $Y$ replaced by $|X|_2^{2p}$, we derive that for every $n\in
\mathbb{N}$
\begin{align*}
    &\mathbb{E}\sup\limits_{s\in[0,t\wedge\tau_n]} |X(s)|_2^p \\
    \leq& |x|_2^p + \sqrt{2|\mu|_{\infty}} pC \mathbb{E} \left[\int_0^{t\wedge\tau_n} |X(s)|_2^{2p}
    ds\right]^{\frac{1}{2}}
     +2p (p-2) |\mu|_{\infty} \mathbb{E} \int_0^{t\wedge\tau_n}
    |X(s)|_2^{p} ds\\
    \leq& |x|_2^p + \epsilon  \sqrt{2|\mu|_{\infty}} pC \mathbb{E}
    \sup\limits_{s\in[0,t\wedge\tau_n]}|X(s)|_2^p
     + C_{\epsilon} \sqrt{2|\mu|_{\infty}} pC \int_0^{t} \mathbb{E} \sup\limits_{r\in[0,s\wedge\tau_n]}|X(r)|_2^p
     ds \\
     &+2p (p-2) |\mu|_{\infty} \int_0^t
    \mathbb{E} \sup\limits_{r\in[0,s\wedge\tau_n]}
    |X(r)|_2^p ds.
\end{align*}

Therefore, similar arguments as at the end of the proof of Theorem
\ref{priesti}
yield Lemma \ref{l4}. \hfill $\square$\\

\end{document}